\documentclass{aims}
\usepackage{amsmath}
\usepackage{paralist}
\usepackage{graphics} 
\usepackage{epsfig} 
\usepackage{graphicx}
\usepackage{epstopdf}
\usepackage[colorlinks=true]{hyperref}
\hypersetup{urlcolor=blue, citecolor=red}

\usepackage[american]{babel}
\usepackage{microtype}
\usepackage{ragged2e}
\allowdisplaybreaks[4]

\usepackage[pagewise]{lineno}

\def\currentvolume{X}
 \def\currentissue{X}
  \def\currentyear{20XX}
   \def\currentmonth{XX}
    \def\ppages{X--XX}

\newtheorem{theorem}{Theorem}[section]
\newtheorem{lemma}[theorem]{Lemma}
\theoremstyle{definition}
\newtheorem{definition}[theorem]{Definition}
\newtheorem{proposition}[theorem]{Proposition}
\newtheorem{assumption}[theorem]{Assumption}
\newtheorem{example}[theorem]{Example}
\newtheorem{remark}{Remark}

\title[SFIDEs with weakly singular kernels] 
      {Stochastic fractional integro-differential equations with weakly singular kernels: well-posedness and Euler--Maruyama approximation}

\author[Xinjie Dai, Aiguo Xiao and Weiping Bu]{}


\email{xinjie@smail.xtu.edu.cn}
\email{xag@xtu.edu.cn}
\email{weipingbu@lsec.cc.ac.cn}

\thanks{This research is supported by the National Natural Science Foundation of China (Nos.\ 12071403, 11671343), the Research Foundation of Education Commission of Hunan Province of China (No.\ 19B565) and the Postgraduate Innovation Fund of Hunan Province in China (No.\ CX20190420).}

\thanks{$^*$ Corresponding author: Aiguo Xiao}

\begin{document}
\maketitle

\centerline{\scshape Xinjie Dai, Aiguo Xiao$^*$ and Weiping Bu}
\medskip
{\footnotesize
 \centerline{School of Mathematics and Computational Science}
 \centerline{Hunan Key Laboratory for Computation and Simulation in Science and Engineering}
 \centerline{Xiangtan University, Hunan 411105, China}
} 

\bigskip

\centerline{(Communicated by the associate editor name)}

\begin{abstract}
This paper considers the initial value problem of general nonlinear stochastic fractional integro-differential equations with weakly singular kernels. Our effort is devoted to establishing some fine estimates to include all the cases of Abel-type singular kernels. Firstly, the existence, uniqueness and continuous dependence on the initial value of the true solution under local Lipschitz condition and linear growth condition are derived in detail. Secondly, the Euler--Maruyama method is developed for solving numerically the equation, and then its strong convergence is proven under the same conditions as the well-posedness. Moreover, we obtain the accurate convergence rate of this method under global Lipschitz condition and linear growth condition. In particular, the Euler--Maruyama method can reach strong first-order superconvergence when $\alpha = 1$. Finally, several numerical tests are reported for verification of the theoretical findings.
\end{abstract}

\section{Introduction}
\label{sec.1}

\indent Integro-differential equations have many influential applications in scientific fields such as biological population \cite{Scudo1971,TeBeest1997}, wave propagation \cite{Lakshmikantham1995} and reactor dynamics \cite{Levin1960}. Because of the increasing development of fractional calculus and the deepening understanding to its non-locality, fractional integro-differential equations appear in electromagnetic wave \cite{Tarasov2009}, population system \cite{Maleki2015,Yuzbasi2013} and other areas. On the other hand, in order to capture ubiquitous noise factors in the actual situation, stochastic integro-differential equations emerge in anomalous diffusion \cite{McKinley2018}, stochastic feedback system \cite{Rao1975} and option pricing \cite{ContTankov2004}. Nowadays, more and more scholars focus on stochastic fractional equations \cite{AnhDoanHuong2019, DoanHuongKloedenVu2020, SonHuong2018, Tuan2020} since it can be applied to explore the memory, hereditary and hidden properties of some noise systems in physics \cite{LiLiu2017}, mathematical finance \cite{Tien2013} and ecological epidemiology \cite{Pedjeu2012,WangXu2016}, etc.

\indent For stochastic fractional integro-differential equations (SFIDEs) with regular kernels, Badr and El-Hoety in \cite{Badr2012} initially discussed the well-posedness of the linear case. Then, block pulse approximation \cite{Asgari2014}, Galerkin methods \cite{Kamrani2015,Kamrani2016,Mohammadi2015}, spectral collocation method \cite{Taheri2017}, operational matrix method \cite{Mirzaee2017} and meshless collocation method \cite{Mirzaee2018} also have been studied. Both the weak singularity of fractional derivative and the low regularity of stochastic noise bring some inevitable difficulties in the concrete analysis. Undoubtedly, it is more involved when the integral kernels are singular, especially for the stochastic integral.

\indent In this paper, motivated by the non-negligible noise source of some practical problems modeled by fractional integro-differential equations with Abel-type singular kernels \cite{ZhangZhu2020}, we consider the following initial value problem of $d$-dimensional nonlinear SFIDEs in It\^{o}'s sense
\begin{equation}\label{eq:1}
\begin{aligned}
D^{\alpha}y(t) =&~ f_0(t,y(t)) + \int_0^t (t-s)^{-\beta_1} f_1(t,s,y(s)) \mathrm{d}s \\
& + \int_0^t (t-s)^{-\beta_2} f_2(t,s,y(s)) \mathrm{d}W(s)
\end{aligned}
\end{equation}
for $t \in \mathcal{J}:=[0,T]$ with $y(0) = y_0 \in \mathbb{R}^d$. Here, $T>0$ is a given real number, $D^{\alpha}$ is the Caputo fractional derivative of order $\alpha\in(0,1]$,~$\beta_1 \in(0,1)$,~$\beta_2 \in(0,\frac{1}{2})$; $W(t)$ denotes an $r$-dimensional standard Wiener process (i.e., Brownian motion) defined on the complete probability space $(\Omega,\mathcal{F},\mathbb{P})$ with a filtration $\{\mathcal{F}_t\}_{t\geq0}$ satisfying the usual conditions (i.e., it is right continuous and $\mathcal{F}_0$ contains all the $\mathbb{P}$-null sets), and the initial value $y_0$ is an $\mathcal{F}_0$-measurable random variable defined on the same probability space such that $\mathbb{E}\big[|y_0|^p\big]<+\infty$, for some integer $p \geq 2$. For the functions $f_i$ ($i=1,2$), we assume that it does not contain the form $(t-s)^{\gamma}$ with $\gamma > 0$. Equation \eqref{eq:1} arises from spatial approximations of stochastic fractional partial differential equations \cite{JinYan2019,LiWang2019,Zou2019} that have been applied to describe random effects on transport of particles in medium with thermal memory \cite{ChenKim2015} and used to optimal controls \cite{Balasubramaniam2017}.

\indent Inspired by the use of Fubini theorem in deterministic fractional integro-differential equations \cite{Aghajani2012,Momani2000}, Dai, Bu and Xiao in \cite{Dai2019} established a connection between SFIDEs with regular kernels (i.e., $\beta_1 = \beta_2 = 0$) and stochastic Volterra integral equations (SVIEs). For the singular SFIDE \eqref{eq:1}, the useful connection can be found at the end of Section \ref{sec.2}. In order to cover all the cases $\alpha\in(0,1]$,~$\beta_1\in(0,1)$ and $\beta_2\in(0,\frac{1}{2})$, we will develop some fine estimates (e.g., Lemmas \ref{le:3.1} and \ref{le:3.2}) and the supreme estimate (see the proof of Lemma \ref{le:3.3}) as well as carefully use H\"{o}lder's inequality and weakly singular Gronwall's inequality. Furthermore, the local Lipschitz condition (see Assumption \ref{as:2.3}) will be used for the nonlinear functions $f_i$ ($i=0,1,2$) to relax the global Lipschitz condition. It is worth noting that It\^{o}'s formula and Doob's martingale inequality, which are frequently employed in the analysis of stochastic differential equations (SDEs), cannot be well applied to the SFIDE \eqref{eq:1} due to their unavailability for singular SVIEs \cite[page 1063]{Wang2008}.

\indent This article is devoted to two main goals:
\begin{enumerate}
\item Under local Lipschitz condition and linear growth condition, the existence and uniqueness as well as continuous dependence on the initial value of the true solution to SFIDE \eqref{eq:1} are derived to fill the well-posedness gap;
\item The strong convergence of the Euler--Maruyama (EM) method is investigated under the same conditions as the well-posedness, and its accurate convergence rate is obtained under global Lipschitz condition and linear growth condition to show its efficiency for solving the SFIDE \eqref{eq:1}.
\end{enumerate}
In particular, when $\alpha = 1$, the considered equation \eqref{eq:1} becomes the integer-order stochastic integro-differential equation, the obtained convergence rate (see Theorem \ref{th.4.6}) actually improves the corresponding result of \cite[the case $H \in (0,\frac{1}{2})$ of Theorem 3.9]{YangYangYao2021}.

\indent The rest of the paper is organized as follows. Some notations and assumptions as well as the connection between SFIDEs and SVIEs will be introduced in Section \ref{sec.2}. Section \ref{sec.3} begins to analyze the well-posedness of the problem \eqref{eq:1}. Section \ref{sec.4} aims to derive the strong convergence properties of the EM method. Several numerical test examples are given in Section \ref{sec.5}. Section \ref{sec.6} affords some brief conclusions.

\section{Mathematical preliminaries}
\label{sec.2}

\indent Throughout this paper, unless otherwise specified, we use the following notations. Let $\mathbb{E}$ denote the expectation corresponding to $\mathbb{P}$. If $A$ is a vector or matrix, then its transpose is denoted by $A^{T}$. Let $|\cdot|$ denote both the Euclidean norm on $\mathbb{R}^d$ and the trace (or Frobenius) norm on $\mathbb{R}^{d\times r}$. That is, if $x\in\mathbb{R}^d$, then $|x|$ is the Euclidean norm; If $A$ is a matrix, then $|A| = \sqrt{\text{trace}(A^{T}A)}$ is its trace norm. If $S$ is a set, then its indicator function is denoted by $1_S$, namely $1_S(x) = 1$ if $x\in S$ and 0 otherwise. For two real numbers $a$ and $b$, we write $a \vee b := \max\{a,b\}$ and $a \wedge b := \min\{a,b\}$. Moreover, the capital letter $C$ (with or without subscripts) will be used to represent a generic positive constant whose value may change when it appears in different places, but it is always independent of the step size $h$. We now impose four mild hypotheses which will be used later for the nonlinear functions $f_i$ ($i=0,1,2$).

\begin{assumption}\label{as:2.1}
There exists a positive constant $L_1$ such that for all $t_1,t_2,s \in \mathcal{J}$ and all $y \in \mathbb{R}^d$, $f_1$ and $f_2$ satisfy the condition:
\begin{align*}
|f_i(t_1,s,y) - f_i(t_2,s,y)| \leq L_1(1+|y|)|t_1 - t_2|, \quad i = 1,2.
\end{align*}
\end{assumption}

\begin{assumption}\label{as:2.2}
There exists a positive constant $L_2$ such that for all $t,s_1,s_2 \in \mathcal{J}$ and all $y \in \mathbb{R}^d$, $f_0$, $f_1$ and $f_2$ satisfy the condition:
\begin{align*}
|f_0(s_1,y) - f_0(s_2,y)| \vee |f_i(t,s_1,y) - f_i(t,s_2,y)| \leq L_2(1+|y|)|s_1 - s_2|, \quad i = 1,2.
\end{align*}
\end{assumption}

\begin{assumption}\label{as:2.3}
For each integer $m\geq1$, there exists a positive constant $K_m$, depending only on $m$, such that for all $t,s\in\mathcal{J}$ and all $y_1,y_2 \in \mathbb{R}^d$ with $|y_1|\vee|y_2| \leq m$, $f_0$, $f_1$ and $f_2$ satisfy the local Lipschitz condition:
\begin{align*}
|f_0(s,y_1) - f_0(s,y_2)| \vee |f_i(t,s,y_1) - f_i(t,s,y_2)| \leq K_m|y_1 - y_2|, \quad i = 1,2.
\end{align*}
\end{assumption}

\begin{assumption}\label{as:2.4}
There exists a positive constant $L$ such that for all $t,s \in \mathcal{J}$ and all $y \in \mathbb{R}^d$, $f_0$, $f_1$ and $f_2$ satisfy the linear growth condition:
\begin{align*}
|f_0(s,y)| \vee |f_i(t,s,y)| \leq L(1 + |y|), \quad i = 1,2.
\end{align*}
\end{assumption}

\begin{remark}
To reflect the generality of our results, we emphasize that the local Lipschitz condition above (i.e., Assumption \ref{as:2.3}) is weaker than the following global Lipschitz condition: there exists a positive constant $K$ such that for all $t,s \in\mathcal{J}$ and all $y_1, y_2 \in \mathbb{R}^d$, $f_0$, $f_1$ and $f_2$ satisfy the inequality
\begin{align}\label{eq:2}
|f_0(s,y_1) - f_0(s,y_2)| \vee |f_i(t,s,y_1) - f_i(t,s,y_2)| \leq K|y_1 - y_2|, \quad i =1, 2.
\end{align}
For example, the function $\cos(y^2)$ satisfies the local Lipschitz condition but does not satisfy the global Lipschitz condition.
\end{remark}

\indent The following preliminary definitions and properties are taken from \cite{AnhDoanHuong2019, Diethelm2010}.

\begin{definition}\label{de:2.5}
Let $\alpha \in (0,1]$, $T \in [0,+\infty)$ and $f: [0,T] \rightarrow \mathbb{R}^d$ be a measurable function such that $\int_{0}^{T} |f(\tau)| \mathrm{d} \tau < +\infty$. The Riemann--Liouville fractional integral operator of order $\alpha$ is defined as
\begin{align*}
I^{\alpha} f(t) = \frac{1}{\Gamma(\alpha)}\int_0^t (t-\tau)^{\alpha-1}f(\tau)\mathrm{d}\tau, \qquad \forall\, t \in [0,T],
\end{align*}
where $\Gamma(\alpha):=\int_0^{+\infty}e^{-t}t^{\alpha-1}\mathrm{d}t$ is the Gamma function.
\end{definition}

\begin{definition}\label{de:2.6}
Let $\alpha \in (0,1]$, $T \in [0,+\infty)$ and $f\in C^1([0,T];\mathbb{R}^d)$. The Caputo derivative of order $\alpha$ can be defined by
$D^{\alpha}f(t) := (I^{1-\alpha} D f)(t)$, where $D = \frac{\mathrm{d}}{\mathrm{d} t}$ is the usual derivative.
\end{definition}

\begin{proposition}\label{pr:2.7}
The Riemann--Liouville fractional integral operator and the Caputo fractional derivative admit the following properties:
\begin{enumerate}
\item $I^{\alpha}(D^{\alpha}f(t)) = f(t) - f(0)$;
\item $D^{\alpha}C = 0$,\quad here $C$ is a constant;
\item $D^{\alpha}(I^{\alpha}f(t)) = f(t)$.
\end{enumerate}
\end{proposition}

\indent As with fractional SDEs (e.g., \cite[Definition 2.2]{AnhDoanHuong2019}, \cite[Page 318]{LiLiu2017}), the SFIDE \eqref{eq:1} is actually rigorously defined by its integral form, namely
\begin{equation}\label{eq:3}
\begin{aligned}
y(t) =&~ y_0 + \frac{1}{\Gamma(\alpha)} \int_0^t (t-\tau)^{\alpha-1} f_0(\tau,y(\tau)) \mathrm{d}\tau \\
& + \frac{1}{\Gamma(\alpha)} \int_0^t (t-\tau)^{\alpha-1} \Big( \int_0^{\tau} (\tau-s)^{-\beta_1} f_1(\tau,s,y(s)) \mathrm{d}s \Big) \mathrm{d}\tau \\
& + \frac{1}{\Gamma(\alpha)} \int_0^t (t-\tau)^{\alpha-1} \Big( \int_0^{\tau} (\tau-s)^{-\beta_2} f_2(\tau,s,y(s)) \mathrm{d}W(s) \Big) \mathrm{d}\tau.
\end{aligned}
\end{equation}
Then, similar to \cite[Definition 2.1 of Chapter 2]{Mao2008}, we proceed to give the definition of the unique solution to SFIDE \eqref{eq:1}.

\begin{definition}\label{de:2.8}
An $\mathbb{R}^d$-valued stochastic process $\{y(t)\}_{t\in\mathcal{J}}$ is called a solution to the SFIDE \eqref{eq:1} if it has the following properties:
\begin{enumerate}
 \item \indent $\{y(t)\}$ is continuous and $\mathcal{F}_t$-adapted;
 \item \indent $f_0 \in \mathcal{L}^1(\mathcal{J}\times\mathbb{R}^d;\mathbb{R}^d)$, $f_1 \in \mathcal{L}^1(\big\{(t,s):0\leq s\leq t\leq T\big\}\times\mathbb{R}^d;\mathbb{R}^d)$ and $f_2 \in \mathcal{L}^2(\big\{(t,s):0\leq s\leq t\leq T\big\}\times\mathbb{R}^d;\mathbb{R}^{d\times r})$;
 \item \indent The equation \eqref{eq:3} almost surely holds for every $t\in\mathcal{J}$.
\end{enumerate}
\indent A solution $\{y(t)\}$ is said to be unique if any other solution $\{\tilde{y}(t)\}$ is indistinguishable from $\{y(t)\}$, that is
\begin{align*}
\mathbb{P}\Big\{ y(t) = \tilde{y}(t) \text{~for all~} t\in\mathcal{J} \Big\} = 1.
\end{align*}
\end{definition}

\indent End this section with the useful connection between SFIDEs and SVIEs.

\begin{theorem}\label{th.2.9}
Let $\alpha\in(0,1]$,~$\beta_1 \in(0,1)$ and $\beta_2 \in(0,\frac{1}{2})$. Then, a continuous and $\mathcal{F}_t$-adapted stochastic process $y(t)$ is a solution of SFIDE \eqref{eq:1}, which is rigorously defined by its integral form \eqref{eq:3}, if and only if it is a solution of the following SVIE
\begin{align}\label{eq:4}
y(t) = y_0 + \int_0^t F_0(t,s,y(s)) \mathrm{d}s + \int_0^t F_1(t,s,y(s)) \mathrm{d}s + \int_0^t F_2(t,s,y(s)) \mathrm{d}W(s),
\end{align}
where
\begin{align*}
& F_0(t,s,y(s)) := \frac{1}{\Gamma(\alpha)} (t-s)^{\alpha-1} f_0(s,y(s)), \\
\ F_i(t,s,y(s)) &:= \frac{1}{\Gamma(\alpha)} \int_s^t (t-\tau)^{\alpha-1} (\tau-s)^{-\beta_i} f_i(\tau,s,y(s)) \mathrm{d}\tau \\
&= \frac{(t-s)^{\alpha-\beta_i}}{\Gamma(\alpha)} \int_0^1 (1-u)^{\alpha-1}u^{-\beta_i} f_i((t-s)u+s,s,y(s)) \mathrm{d}u, \quad i=1,2.
\end{align*}
\end{theorem}

\begin{proof}
Let $y(t)$ be a solution of SFIDE \eqref{eq:1}, so Eq.\ \eqref{eq:3} almost surely holds and $f_2 \in \mathcal{L}^2(\big\{(t,s):0\leq s\leq t\leq T\big\}\times\mathbb{R}^d;\mathbb{R}^{d\times r})$. Then, it follows from \cite[Theorem 2.1]{Diethelm2010} that
\begin{align*}
\int_0^t (t-\tau)^{\alpha-1} \Big( \int_0^{\tau} (\tau-s)^{-2\beta_2} \mathbb{E}\Big[ \big| f_2(\tau,s,y(s)) \big|^2\Big] \mathrm{d}s \Big)^{\frac{1}{2}} \mathrm{d}\tau < +\infty, \quad \forall\, t\in\mathcal{J},
\end{align*}
which shows that the sufficient condition of stochastic Fubini theorem \cite[Theorem 4.33 or 5.10]{DaPrato2014} is satisfied. Thus, applying Fubini's theorem to Eq.\ \eqref{eq:3} yields
\begin{align*}
y(t) =&~y_0 + \frac{1}{\Gamma(\alpha)} \int_0^t (t-s)^{\alpha-1} f_0(s,y(s)) \mathrm{d}s \\
& + \frac{1}{\Gamma(\alpha)} \int_0^t \Big( \int_s^t (t-\tau)^{\alpha-1} (\tau-s)^{-\beta_1} f_1(\tau,s,y(s)) \mathrm{d}\tau \Big) \mathrm{d}s \\
& + \frac{1}{\Gamma(\alpha)} \int_0^t \Big( \int_s^t (t-\tau)^{\alpha-1} (\tau-s)^{-\beta_2} f_2(\tau,s,y(s)) \mathrm{d}\tau \Big) \mathrm{d}W(s),
\end{align*}
which implies that $y(t)$ is also a solution to the SVIE \eqref{eq:4}.

In reverse, let $y(t)$ be a solution of SVIE \eqref{eq:4}. Then, $y(t)$ is a solution to Eq.\ \eqref{eq:3} by Fubini's theorem, which implies that $y(t)$ is also a solution of SFIDE \eqref{eq:1}.
\end{proof}

\section{Well-posedness of SFIDE \eqref{eq:1}}
\label{sec.3}

\indent Through the preparation of the previous section, we will investigate the existence, uniqueness and continuous dependence on the initial value of the exact solution to SFIDE \eqref{eq:1}.

\subsection{Existence and uniqueness theorem}

\indent To facilitate the proof of the existence result, we first introduce the EM approximation. For every integer $N\geq1$, the EM approximation \cite{Kloeden1992,Mao2008} can be shown as
\begin{equation}\label{eq:5}
\begin{aligned}
y^N(t) =&~y_0 + \int_0^t F_0(t,s,\hat{y}^N(s)) \mathrm{d}s + \int_0^t F_1(t,s,\hat{y}^N(s)) \mathrm{d}s \\
& + \int_0^t F_2(t,s,\hat{y}^N(s)) \mathrm{d}W(s),
\end{aligned}
\end{equation}
where the simple step process $\hat{y}^N(t) = \sum_{n=0}^N y^N(t_n) 1_{[t_n,t_{n+1})}(t)$ and the mesh points $t_n = nh$ ($n=0,1,\cdots,N$) with $h=\frac{T}{N}$.

\indent Next, the following four useful lemmas will be prepared.

\begin{lemma}\label{le:3.1}
Let $0\leq t' < t \leq T$ and $a,b\in(0,1]$. Then,
\begin{align*}
\int_0^{t'} \big|(t - s)^{a-b} - (t' - s)^{a-b}\big| \mathrm{d}s \leq C (t-t')^{1 \wedge (1+a-b)},
\end{align*}
where the positive constant $C$ only depends on $T,a,b$.
\end{lemma}

\begin{proof}
Setting $a,b\in(0,1]$ implies $a-b\in(-1,1)$. Firstly, for $a-b\in(0,1)$, one can derive
\begin{align*}
&\quad \int_0^{t'} \big|(t - s)^{a-b} - (t' - s)^{a-b}\big| \mathrm{d}s = \int_0^{t'} (t - s)^{a-b} - (t' - s)^{a-b} \mathrm{d}s \\
&= \frac{1}{1+a-b} \left( t^{1+a-b} - (t-t')^{1+a-b} - t'^{1+a-b} \right) \leq \frac{1}{1+a-b} \left( t^{1+a-b} - t'^{1+a-b} \right) \\
&= \int_{t'}^{t} s^{a-b} \mathrm{d}s \leq \int_{t'}^{t} t^{a-b} \mathrm{d}s \leq T^{a-b} (t-t').
\end{align*}
Secondly, for $a-b \in (-1,0]$, one can read
\begin{align*}
&\quad\ \int_0^{t'} \big|(t - s)^{a-b} - (t' - s)^{a-b}\big| \mathrm{d}s = \int_0^{t'} (t' - s)^{a-b} - (t - s)^{a-b} \mathrm{d}s \\
&= \frac{1}{1+a-b} \left( t'^{1+a-b} - t^{1+a-b} + (t-t')^{1+a-b} \right) \leq \frac{1}{1+a-b} (t-t')^{1+a-b}.
\end{align*}
The proof is complete.
\end{proof}

\begin{lemma}\label{le:3.2}
Let $0\leq t' < t \leq T$, $a\in(0,1]$ and $b\in(0,\frac{1}{2})$. Then, $a-b\in(-\frac{1}{2},1)$ and
\begin{align*}
\int_0^{t'} \big|(t - s)^{a-b} - (t' - s)^{a-b}\big|^2 \mathrm{d}s \leq
\begin{cases}
C (t-t')^2, &\text{\quad if } a-b \in (\frac{1}{2},1), \\
C (t-t')^{2-\epsilon}, &\text{\quad if } a-b = \frac{1}{2}, \\
C (t-t')^{1+2(a-b)}, &\text{\quad if } a-b \in (-\frac{1}{2},\frac{1}{2}),
\end{cases}
\end{align*}
where the positive constant $C$ only depends on $T,a,b,\epsilon$, and $\epsilon \in (0,1)$.
\end{lemma}

\begin{proof}
If $a-b \in (\frac{1}{2},1)$, then it follows from $a-b-1<0$ and $2(a-b)-1 > 0$ that
\begin{align*}
&\quad\ \int_0^{t'} \big|(t - s)^{a-b} - (t' - s)^{a-b}\big|^2 \mathrm{d}s = \int_0^{t'} \big| (a-b) \int_{t'}^t (\tau - s)^{a-b-1} \mathrm{d}\tau \big|^2 \mathrm{d}s \\
& \leq \int_0^{t'} \big| (a-b) \int_{t'}^t (t' - s)^{a-b-1} \mathrm{d}\tau \big|^2 \mathrm{d}s = (a-b)^2 (t-t')^2 \int_0^{t'} (t' - s)^{2(a-b-1)} \mathrm{d}s \\
&\leq \frac{(a-b)^2 T^{2(a-b)-1}}{2(a-b)-1} (t-t')^2.
\end{align*}
If $a-b \in (0,\frac{1}{2})$, then Cauchy--Schwarz's inequality indicates
\begin{align*}
&\quad\ \int_0^{t'} \big|(t - s)^{a-b} - (t' - s)^{a-b}\big|^2 \mathrm{d}s = \int_0^{t'} \big| (a-b) \int_{t'}^t (\tau - s)^{a-b-1} \mathrm{d}\tau \big|^2 \mathrm{d}s \\
& \leq (a-b)^2 (t-t') \int_0^{t'} \int_{t'}^t (\tau - s)^{2(a-b-1)} \mathrm{d}\tau \mathrm{d}s \\
&= \frac{(a-b)^2}{1-2(a-b)} (t-t') \int_0^{t'} (t' - s)^{2(a-b)-1} - (t - s)^{2(a-b)-1} \mathrm{d}s \\
&= \frac{(a-b)}{2-4(a-b)} (t-t') \left( t'^{2(a-b)} - t^{2(a-b)} + (t - t')^{2(a-b)} \right) \\
&\leq \frac{(a-b)}{2-4(a-b)} (t - t')^{1+2(a-b)}.
\end{align*}
If $a-b = \frac{1}{2}$, then for any $\epsilon \in (0,1)$, Cauchy--Schwarz's inequality shows
\begin{align*}
&\quad\ \int_0^{t'} \big|(t - s)^{a-b} - (t' - s)^{a-b}\big|^2 \mathrm{d}s = \int_0^{t'} \big| \frac{1}{2} \int_{t'}^t (\tau - s)^{-\frac{1}{2}} \mathrm{d}\tau \big|^2 \mathrm{d}s \\
&\leq \frac{1}{4} \int_0^{t'} \int_{t'}^t (\tau - s)^{-\epsilon} \mathrm{d}\tau \int_{t'}^t (\tau - s)^{-1+\epsilon} \mathrm{d}\tau \mathrm{d}s \\
&= \frac{1}{4(1-\epsilon)} \int_0^{t'} \left( (t-s)^{1-\epsilon} - (t'-s)^{1-\epsilon} \right) \int_{t'}^t (\tau - s)^{-1+\epsilon} \mathrm{d}\tau \mathrm{d}s \\
&\leq \frac{1}{4(1-\epsilon)} (t-t')^{1-\epsilon} \int_0^{t'} \int_{t'}^t (t' - s)^{-1+\epsilon} \mathrm{d}\tau \mathrm{d}s \\
&= \frac{1}{4(1-\epsilon)} (t-t')^{2-\epsilon} \int_0^{t'} (t' - s)^{-1+\epsilon} \mathrm{d}s \leq C (t-t')^{2-\epsilon}.
\end{align*}
If $a-b \in (-\frac{1}{2},0]$, then it follows from Lemma \ref{le:3.1} that
\begin{align*}
&\quad\ \int_0^{t'} \big|(t - s)^{a-b} - (t' - s)^{a-b}\big|^2 \mathrm{d}s \\
&=\int_0^{t'}(t - s)^{2(a-b)} - 2 (t - s)^{a-b} (t' - s)^{a-b} + (t' - s)^{2(a-b)} \mathrm{d}s \\
&\leq \int_0^{t'}(t - s)^{2(a-b)} - 2 (t - s)^{a-b} (t - s)^{a-b} + (t' - s)^{2(a-b)} \mathrm{d}s \\
&= \int_0^{t'} (t' - s)^{2(a-b)} - (t - s)^{2(a-b)} \mathrm{d}s \leq \frac{1}{1+2(a-b)} (t-t')^{1+2(a-b)}.
\end{align*}
Hereto, Lemma \ref{le:3.2} has been proven.
\end{proof}

\begin{lemma}\label{le:3.3}
If the Assumption \ref{as:2.4} holds, then there exists a positive constant $C$ independent of $N$ such that for any integer $p \geq 2$,
\begin{align*}
\mathbb{E}\Big[|y^N(t)|^p\Big] \leq C \quad \text{and} \quad \mathbb{E}\Big[|\hat{y}^N(t)|^p\Big] \leq C, \qquad \forall\, t\in\mathcal{J}.
\end{align*}
\end{lemma}

\begin{proof}
We first prove the case of $p>2$. For every integer $k\geq1$, we define the stopping time
\begin{align*}
\rho_k^N = T \wedge \inf\{ t\in\mathcal{J}: |y^N(t)|\geq k \},
\end{align*}
where $\rho_k^N \uparrow T$ almost surely as $k\rightarrow +\infty$. For convenience, set $y_k^N(t) = y^N(t\wedge\rho_k^N)$ and $\hat{y}_k^N(t) = \hat{y}^N(t\wedge\rho_k^N)$ for all $t\in\mathcal{J}$. With the help of H\"{o}lder's inequality, BDG's inequality as well as Assumption \ref{as:2.4} and $\mathbb{E}\big[|y_0|^p\big]<+\infty$, one can derive from \eqref{eq:5} that there is a positive constant $q\in(0, p \alpha)$, which only depends on $p>2$ and $\alpha\in(0,1]$, such that
\begin{align*}
& \mathbb{E}\Big[|y_k^N(t)|^p\Big] \leq \frac{4^{p-1}}{\Gamma^p(\alpha)} \Bigg\{ \Gamma^p(\alpha)\mathbb{E}\Big[|y_0|^p\Big] + \mathbb{E}\Big[ \big| \int_0^{t\wedge\rho_k^N} (t\wedge\rho_k^N - s)^{\alpha-1} f_0(s,\hat{y}_k^N(s)) \mathrm{d}s \big|^p \Big]\\
&\quad + B^p(\alpha,1-\beta_1) \mathbb{E} \Big[\big| \int_0^{t\wedge\rho_k^N} (t\wedge\rho_k^N-s)^{\alpha-\beta_1} \sup_{s\leq u\leq (t\wedge\rho_k^N)} |f_1(u,s,\hat{y}_k^N(s))| \mathrm{d}s \big|^p \Big] \\
&\quad + B^p(\alpha,1-\beta_2) \mathbb{E} \Big[ \Big| \int_0^{t\wedge\rho_k^N} (t\wedge\rho_k^N -s)^{2(\alpha-\beta_2)} \sup_{s\leq u\leq (t\wedge\rho_k^N)} |f_2(u,s,\hat{y}_k^N(s))|^2 \mathrm{d}s \Big|^{\frac{p}{2}} \Big] \Bigg\} \\
&\leq C_1 \Bigg\{ 1 + 2\left( \int_0^{t\wedge\rho_k^N} (t\wedge\rho_k^N - s)^{\frac{p\alpha-q}{p-1}-1} \mathrm{d}s \right)^{p-1} \cdot \\
&\qquad\qquad\qquad \int_0^{t\wedge\rho_k^N} (t\wedge\rho_k^N - s)^{q-1} \left( 1 + \mathbb{E}\Big[|\hat{y}_k^N(s)|^p\Big] \right) \mathrm{d}s \\
&\qquad\ \ + \left( \int_0^{t\wedge\rho_k^N} (t\wedge\rho_k^N - s)^{\frac{2(p\alpha-q)}{p-2}-1} \mathrm{d}s \right)^{\frac{p-2}{2}} \cdot \\
&\qquad\qquad\qquad \int_0^{t\wedge\rho_k^N} (t\wedge\rho_k^N - s)^{q-1} \left( 1 + \mathbb{E}\Big[|\hat{y}_k^N(s)|^p\Big] \right) \mathrm{d}s \Bigg\} \\
&\leq C_2 \left( 1 + \int_0^{t\wedge\rho_k^N} (t\wedge\rho_k^N - s)^{q-1} \mathbb{E}\Big[|\hat{y}_k^N(s)|^p\Big] \mathrm{d}s \right),
\end{align*}
where $B(a,b) := \int_0^1 (1-u)^{a-1}u^{b-1} \mathrm{d}u$ with $a,b>0$ is the Beta function, and the positive constants $C_1$ and $C_2$ are independent of $k$ and $N$. Taking the supremum on both sides shows
\begin{align*}
&\quad\ \sup_{0\leq r\leq t}\mathbb{E}\Big[|y_k^N(r)|^p\Big] \\
&\leq C_2 \left( 1 + \sup_{0\leq r\leq t} \int_0^{r\wedge\rho_k^N} (r\wedge\rho_k^N - s)^{q-1} \sup_{0\leq \xi\leq s} \mathbb{E}\Big[|y_k^N(\xi)|^p\Big] \mathrm{d}s \right) \\
&= C_2 \left( 1 + \sup_{0\leq r\leq t} (r\wedge\rho_k^N)^q \int_0^1 (1-u)^{q-1} \sup_{0\leq \xi\leq (r\wedge\rho_k^N)u} \mathbb{E}\Big[|y_k^N(\xi)|^p\Big] \mathrm{d}u \right) \\
&\qquad\qquad\qquad\qquad\qquad\qquad\qquad\qquad\qquad\qquad\qquad \text{(where $u := \frac{s}{r\wedge\rho_k^N}$)} \\
&\leq C_2 \left( 1 + t^q \int_0^1 (1-u)^{q-1} \sup_{0\leq \xi\leq tu} \mathbb{E}\Big[|y_k^N(\xi)|^p\Big] \mathrm{d}u \right) \\
&= C_2 \left( 1 + \int_0^t (t-s)^{q-1} \sup_{0\leq \xi\leq s} \mathbb{E}\Big[|y_k^N(\xi)|^p\Big] \mathrm{d}s \right) \quad \text{(where $s := tu$)},
\end{align*}
which with the application of weakly singular Gronwall's inequality \cite[Corollary 2]{YeGao2007} yields
\begin{align*}
\mathbb{E}\Big[|y_k^N(t)|^p\Big] \leq C, \qquad \forall\, t\in\mathcal{J}.
\end{align*}
Letting $k\rightarrow+\infty$ and using Fatou's lemma to indicate
\begin{align*}
\mathbb{E}\Big[|y^N(t)|^p\Big] \leq C, \qquad \forall\, t\in\mathcal{J},
\end{align*}
which also implies
\begin{align*}
\mathbb{E}\Big[|\hat{y}^N(t)|^p\Big] \leq C, \qquad \forall\, t\in\mathcal{J}.
\end{align*}

Secondly, based on the above proof for the case $p>2$, it is easy to obtain the same conclusions for the case $p=2$, where only H\"{o}lder's inequality is replaced by Cauchy--Schwarz's inequality. The proof is completed.
\end{proof}

\begin{lemma}\label{le:3.4}
Let $0\leq t' < t \leq T$. If the Assumptions \ref{as:2.1} and \ref{as:2.4} hold, then there exists a positive constant $C$ independent of $N$ such that for any integer $p \geq 2$,
\begin{align*}
\mathbb{E}\Big[|y^N(t) - y^N(t')|^p\Big] \leq C (t-t')^{\alpha p}.
\end{align*}
\end{lemma}

\begin{proof}
It follows from \eqref{eq:5} that
\begin{align*}
& \mathbb{E}\Big[|y^N(t) - y^N(t')|^p\Big] \leq 3^{p-1} \bigg\{ \mathbb{E}\Big[ \big| \int_{0}^{t} F_0(t,s,\hat{y}^N(s)) \mathrm{d}s - \int_{0}^{t'} F_0(t',s,\hat{y}^N(s)) \mathrm{d}s \big|^p \Big] \\
&\qquad\qquad\qquad + \mathbb{E}\Big[ \big| \int_{0}^{t} F_1(t,s,\hat{y}^N(s)) \mathrm{d}s - \int_{0}^{t'} F_1(t',s,\hat{y}^N(s)) \mathrm{d}s \big|^p \Big] \\
&\qquad\qquad\qquad + \mathbb{E}\Big[ \big| \int_{0}^{t} F_2(t,s,\hat{y}^N(s)) \mathrm{d}W(s) - \int_{0}^{t'} F_2(t',s,\hat{y}^N(s)) \mathrm{d}W(s) \big|^p \Big] \bigg\} \\
&\qquad\qquad\quad =: 3^{p-1} \{ \mathcal{K}_0 + \mathcal{K}_1 + \mathcal{K}_2 \}.
\end{align*}
Using H\"{o}lder's inequality and Assumptions \ref{as:2.1}, \ref{as:2.4} as well as Lemmas \ref{le:3.1} and \ref{le:3.3} shows
\begin{align*}
\mathcal{K}_0 &\leq 2^{p-1} \bigg\{ \mathbb{E}\Big[ \big| \int_{0}^{t'} \big( (t-s)^{\alpha-1} - (t'-s)^{\alpha-1} \big) f_0(s,\hat{y}^N(s)) \mathrm{d}s \big|^p \Big] \\
&\qquad\qquad + \mathbb{E}\Big[ \big| \int_{t'}^{t} (t-s)^{\alpha-1}f_0(s,\hat{y}^N(s)) \mathrm{d}s \big|^p \Big] \bigg\} \\
&\leq C \bigg\{ \Big( \int_0^{t'} |(t-s)^{\alpha-1} - (t'-s)^{\alpha-1}| \mathrm{d} s \Big)^{p-1} \cdot \\
&\qquad\qquad\qquad \int_0^{t'} |(t-s)^{\alpha-1} - (t'-s)^{\alpha-1}| \left(1+\mathbb{E}\Big[|\hat{y}^N(s)|^p\Big]\right) \mathrm{d} s \\
&\qquad\quad + \Big( \int_{t'}^t (t-s)^{\alpha-1} \mathrm{d}s \Big)^{p-1} \cdot \int_{t'}^t (t-s)^{\alpha-1} \left(1+\mathbb{E}\Big[|\hat{y}^N(s)|^p\Big]\right) \mathrm{d} s \bigg\} \\
&\leq C (t-t')^{\alpha p},
\end{align*}
and
\begin{align*}
\mathcal{K}_1 \leq C (t-t')^{(1 \wedge (1 + \alpha - \beta_1)) p}.
\end{align*}
Applying H\"{o}lder's inequality, BDG inequality and Assumptions \ref{as:2.1}, \ref{as:2.4} as well as Lemmas \ref{le:3.2} and \ref{le:3.3} displays
\begin{align*}
\mathcal{K}_2 &\leq C \bigg\{ \mathbb{E}\Big[ \big| \int_0^{t'} \big( (t-s)^{\alpha-\beta_2} - (t'-s)^{\alpha-\beta_2} \big)^2 \sup_{s\leq u\leq t}\big|f_2(u,s,\hat{y}^N(s))\big|^2 \mathrm{d}s \big|^{\frac{p}{2}} \Big] \\
&\qquad\quad + \mathbb{E}\Big[ \big| \int_0^{t'} (t'-s)^{2(\alpha-\beta_2)} \sup_{0\leq u\leq 1} \big|f_2((t-s)u+s,s,\hat{y}^N(s)) - \\
&\qquad\qquad\qquad\qquad\qquad\qquad\qquad\qquad\ \ f_2((t'-s)u+s,s,\hat{y}^N(s)) \big|^2 \mathrm{d}s \big|^{\frac{p}{2}} \Big] \\
&\qquad\quad + \mathbb{E}\Big[ \big| \int_{t'}^t (t-s)^{2(\alpha-\beta_2)} \sup_{s\leq u\leq t}\big|f_2(u,s,\hat{y}^N(s))\big|^2 \mathrm{d}s \big|^{\frac{p}{2}} \bigg\} \\
&\leq C \bigg\{ \big| \int_0^{t'} \big( (t-s)^{\alpha-\beta_2} - (t'-s)^{\alpha-\beta_2} \big)^2 \mathrm{d}s \big|^{\frac{p-2}{2}} \cdot \\
&\qquad\qquad \int_0^{t'} \big( (t-s)^{\alpha-\beta_2} - (t'-s)^{\alpha-\beta_2} \big)^2 \left(1+\mathbb{E}\Big[|\hat{y}^N(s)|^p\Big]\right) \mathrm{d}s \\
&\quad + (t-t')^2 \big| \int_0^{t'} (t'-s)^{2(\alpha-\beta_2)} \mathrm{d}s \big|^{\frac{p-2}{2}} \int_0^{t'} (t'-s)^{2(\alpha-\beta_2)} \left(1+\mathbb{E}\Big[|\hat{y}^N(s)|^p\Big]\right) \mathrm{d}s \\
&\quad + \big| \int_{t'}^t (t-s)^{2(\alpha-\beta_2)} \mathrm{d}s \big|^{\frac{p-2}{2}} \int_{t'}^t (t-s)^{2(\alpha-\beta_2)} \left(1+\mathbb{E}\Big[|\hat{y}^N(s)|^p\Big]\right) \mathrm{d}s \bigg\} \\
&\leq
\begin{cases}
C(t-t')^{(2-\epsilon)\frac{p}{2}}, & \mbox{if } \alpha - \beta_2 = \frac{1}{2}, \\
C (t-t')^{(2 \wedge (1 + 2\alpha - 2\beta_2)) \frac{p}{2}}, & \mbox{otherwise}.
\end{cases}
\end{align*}
If $\alpha - \beta_2 = \frac{1}{2}$, then $2-2\alpha \in (0,1)$ since $\beta_2 \in (0,\frac{1}{2})$. Taking $\epsilon = 2-2\alpha$ yields
\begin{align*}
\mathcal{K}_2 \leq C (t-t')^{\alpha p}.
\end{align*}
Therefore,
\begin{align*}
\mathbb{E}\Big[|y^N(t) - y^N(t')|^p\Big] \leq 3^{p-1} \big( \mathcal{K}_0 + \mathcal{K}_1 + \mathcal{K}_2 \big) \leq C (t-t')^{\alpha p},
\end{align*}
which completes the proof.
\end{proof}

\begin{theorem}\label{th:3.5}
Under the Assumptions \ref{as:2.1}, \ref{as:2.3} and \ref{as:2.4}, the SFIDE \eqref{eq:1} has a unique solution $y(t)$. Moreover, for any positive integer $p\geq2$,
\begin{align*}
\mathbb{E}\Big[|y(t)|^p\Big] < +\infty, \qquad \forall\, t\in\mathcal{J}.
\end{align*}
\end{theorem}

\begin{proof}
For the sake of simplicity, we only prove the case of the global Lipschitz condition (i.e., Assumption \ref{as:2.3} is replaced with the condition \eqref{eq:2}). In fact, based on the proof under the global Lipschitz condition, by mean of the truncation functions
\begin{equation*}
F_i^m(t,s,y):=\left\{
\begin{aligned}
&F_i(t,s,y) \qquad &\text{if~} |y|\leq m, \nonumber\\
&F_i(t,s,\frac{my}{|y|}) \qquad &\text{if~} |y| > m,
\end{aligned}
\right.
\qquad
i=0,1,2
\end{equation*}
for each integer $m\geq1$, one can easily read the desired result for the case of the local Lipschitz condition by the similar proof procedure of \cite[Theorem 3.4 of Chapter 2]{Mao2008}.

\indent\underline{Uniqueness}. Let $y(t)$ and $\tilde{y}(t)$ be two solutions of the SFIDE \eqref{eq:1} on the same probability space with $y(0) = \tilde{y}(0)$. Then, Theorem \ref{th.2.9} shows that $y(t)$ and $\tilde{y}(t)$ are also two solutions to the SVIE \eqref{eq:4}. Using H\"{o}lder's inequality, It\^{o} isometry as well as the Lipschitz condition \eqref{eq:2}, one can read from \eqref{eq:4} that
\begin{align}\label{eq:6}
\mathbb{E}\Big[|y(t) - \tilde{y}(t)|^2\Big] &\leq C \mathbb{E} \Bigg\{ \Big|\int_0^t (t-s)^{\alpha-1} \Big( f_0(s,y(s)) - f_0(s,\tilde{y}(s)) \Big) \mathrm{d}s \Big|^2 \nonumber\\
&\qquad\ + \Big|\int_0^t (t-s)^{\alpha-\beta_1} \sup_{s\leq u\leq t} \big|f_1(u,s,y(s)) - f_1(u,s,\tilde{y}(s))\big| \mathrm{d}s \Big|^2 \nonumber\\
&\qquad\ + \int_0^t (t-s)^{2(\alpha-\beta_2)} \sup_{s\leq u\leq t} \big|f_2(u,s,y(s)) - f_2(u,s,\tilde{y}(s))\big|^2 \mathrm{d}s \Bigg\} \nonumber\\
&\leq C \int_0^t (t-s)^{\alpha-1} \mathbb{E}\Big[|y(s) - \tilde{y}(s)|^2\Big] \mathrm{d}s.
\end{align}
Then, weakly singular Gronwall's inequality \cite[Corollary 2]{YeGao2007} yields
\begin{align*}
\mathbb{E}\Big[|y(t) - \tilde{y}(t)|^2\Big] = 0, \qquad \forall\, t\in\mathcal{J},
\end{align*}
which indicates
\begin{align*}
\mathbb{P} \Big\{ |y(t) - \tilde{y}(t)| = 0,~\forall\, t \in \mathbb{Q} \cap \mathcal{J}\Big\} = 1,
\end{align*}
where $\mathbb{Q}$ represents the set of all rational numbers. It follows from the continuity of $|y(t) - \tilde{y}(t)|$ with respect to $t$ that
\begin{align*}
\mathbb{P} \Big\{ |y(t) - \tilde{y}(t)| = 0,~\forall\, t \in \mathcal{J}\Big\} = 1.
\end{align*}
The uniqueness has been proven.

\indent\underline{Existence}. Let $M\geq N\geq1$. Similar to the derivation of the estimate \eqref{eq:6}, one can claim from \eqref{eq:5} that for any $p \geq 2$,
\begin{align*}
& \mathbb{E}\Big[|y^M(t) - y^N(t)|^p\Big] \leq C\int_0^t (t-s)^{\alpha-1} \mathbb{E}\Big[|\hat{y}^M(s) - \hat{y}^N(s)|^p\Big] \mathrm{d}s,
\end{align*}
where the positive constant $C$ is independent of $M$ and $N$. Next, according to the triangle inequality with Lemma \ref{le:3.4} and weakly singular Gronwall's inequality, one concludes that $\{y^N(t)\}$ is a Cauchy sequence and has a limit $y(t)$ in $\mathcal{L}^p\big(\Omega;\mathbb{R}^d\big)$. Obviously, $y(t)$ is $\mathcal{F}_t$-adapted. Moreover, it follows from Lemma \ref{le:3.4} and Fatou's lemma that
\begin{align*}
\mathbb{E}\Big[|y(t) - y(t')|^p\Big] \leq C (t-t')^{\alpha p}.
\end{align*}
By Kolmogorov's continuity criterion, the process $y(t)$ has a continuous version. Letting $N\rightarrow+\infty$ in \eqref{eq:5} indicates that the continuous version is a solution to the SVIE \eqref{eq:4} on $[0,T]$. Then, Theorem \ref{th.2.9} implies that the continuous version is also a solution to SFIDE \eqref{eq:1}. Finally, by Lemma \ref{le:3.3}, letting $N\rightarrow+\infty$ yields the conclusion
\begin{align*}
\mathbb{E}\Big[|y(t)|^p\Big] < +\infty, \qquad \forall\, t\in\mathcal{J}.
\end{align*}
Hereto, Theorem \ref{th:3.5} has been proven.
\end{proof}

\subsection{Continuous dependence of solutions on the initial value}

\begin{definition}\label{de:3.6}
The solution of the SFIDE \eqref{eq:1} is said to be continuous depending on the initial value in mean square sense, if for any $\varepsilon$, there exists a positive number $\zeta$ such that
\begin{align*}
\mathbb{E}\Big[|y(t) - z(t)|^2\Big] < \varepsilon, \qquad \forall\, t \in \mathcal{J},
\end{align*}
provided that $\mathbb{E}\big[|y_0 - z_0|^2\big] < \zeta$, where $z(t)$ is any other solution to \eqref{eq:1} with the initial value $z_0 \in \mathbb{R}^d$.
\end{definition}

\begin{theorem}\label{th:3.7}
Under the Assumptions \ref{as:2.1}, \ref{as:2.3} and \ref{as:2.4}, the solution of SFIDE \eqref{eq:1} continuously depends on the initial value in mean square sense.
\end{theorem}

\begin{proof}
Let $y(t)$ and $z(t)$ be two solutions to \eqref{eq:1} with different initial values $y_0$ and $z_0$, respectively. Then, $y(t)$ and $z(t)$ are also two solutions to \eqref{eq:4}. For simplicity, let $e(t) = y(t) - z(t)$,
\begin{align}\label{eq:7}
\rho_m = \inf \big\{t\geq0: |y(t)|\geq m \big\}, \quad \nu_m = \inf \big\{t\geq0: |z(t)|\geq m \big\}, \quad \sigma_m = \rho_m \wedge\nu_m
\end{align}
for each integer $m\geq1$. Recall Young's inequality: for $u,v>1$ with $u^{-1}+v^{-1} = 1$,
\begin{align*}
ab \leq \frac{\delta}{u}a^u + \frac{1}{v\delta^{\frac{v}{u}}}b^v, \qquad \forall\, a,b,\delta>0.
\end{align*}
Then it holds that for any $\delta > 0$,
\begin{equation}\label{eq:8}
\begin{aligned}
&\quad\ \mathbb{E} \Big[|e(t)|^2\Big] = \mathbb{E}\Big[|e(t)|^2 1_{\{\rho_m>T,~\nu_m>T\}}\Big] + \mathbb{E}\Big[|e(t)|^2 1_{\{\rho_m\leq T \text{~or~} \nu_m\leq T\}}\Big] \\
&\leq \mathbb{E} \Big[ |e(t\wedge\sigma_m)|^2 1_{\{\sigma_m>T\}} \Big] + \frac{2\delta}{p} \mathbb{E} \Big[|e(t)|^p\Big] + \frac{p-2}{p\delta^{\frac{2}{p-2}}}\mathbb{P}\Big( \rho_m\leq T \text{~or~} \nu_m\leq T \Big).
\end{aligned}
\end{equation}
On the one hand, Theorem \ref{th:3.5} gives
\begin{align*}
\mathbb{E} \Big[|e(t)|^p\Big] \leq 2^{p-1}\mathbb{E}\Big[|y(t)|^p + |z(t)|^p\Big] \leq 2^p M,
\end{align*}
here and in the rest of the proof, $M$ stands for a positive constant independent of $m$ and $\delta$. On the other hand,
\begin{align*}
\mathbb{P}\Big( \rho_m\leq T \text{~or~} \nu_m\leq T \Big) \leq \mathbb{P}\Big( \rho_m\leq T\Big) + \mathbb{P}\Big( \nu_m\leq T \Big) \leq \frac{2M}{m^p},
\end{align*}
where we employed the estimate
\begin{align*}
\mathbb{P}\Big( \rho_m\leq T\Big) = \mathbb{E}\Big[1_{\{\rho_m\leq T\}} \frac{|y(\rho_m)|^p}{m^p}\Big] \leq \frac{1}{m^p} \mathbb{E} \Big[|y(\rho_m \wedge T)|^p\Big] \leq \frac{M}{m^p}
\end{align*}
and a same estimate for $\mathbb{P}\big( \nu_m\leq T\big)$. As thus, the inequality \eqref{eq:8} gives
\begin{align}\label{eq:9}
\mathbb{E} \Big[|e(t)|^2\Big] \leq \mathbb{E} \Big[|e(t\wedge\sigma_m)|^2\Big] + \frac{2^{p+1}\delta M}{p} + \frac{(p-2)2M}{p\delta^{\frac{2}{p-2}}m^p}.
\end{align}
Using a similar derivation of the inequality \eqref{eq:6} yields
\begin{align*}
&\quad\ \mathbb{E}\Big[|e(t\wedge\sigma_m)|^2\Big] \\
&\leq CK_m^2 \left( \mathbb{E}\Big[|y_0 - z_0|^2\Big] + \int_0^{t\wedge\sigma_m} ({t\wedge\sigma_m}-s)^{\alpha-1} \mathbb{E}\Big[|e(s\wedge\sigma_m)|^2\Big] \mathrm{d}s \right).
\end{align*}
Applying weakly singular Gronwall's inequality and then inserting the obtained result into \eqref{eq:9} lead to
\begin{align*}
\mathbb{E}\Big[|e(t)|^2\Big] \leq C_m \mathbb{E}\Big[|y_0 - z_0|^2\Big] + \frac{2^{p+1}\delta M}{p} + \frac{(p-2)2M}{p\delta^{\frac{2}{p-2}}m^p},
\end{align*}
where the positive constant $C_m$ depends on $m$, but not on $\delta$. It is observed that $\delta$ and $m$ can be can chosen such that
\begin{align*}
\frac{2^{p+1}\delta M}{p} < \frac{\varepsilon}{3}, \qquad \frac{(p-2)2M}{p\delta^{\frac{2}{p-2}}m^p} < \frac{\varepsilon}{3},
\end{align*}
and there exists a positive number $\zeta>0$ such that $C_m \mathbb{E}\big[|y_0 - z_0|^2\big] < C_m \zeta < \frac{\varepsilon}{3}$. Consequently,
\begin{align*}
\mathbb{E}\Big[|e(t)|^2\Big] < \varepsilon, \qquad \forall\, t \in \mathcal{J}.
\end{align*}
The proof is complete.
\end{proof}

\section{Strong convergence and convergence rate of the EM method}
\label{sec.4}

\indent Since it is difficult to obtain the closed-form solution to the SFIDE \eqref{eq:1}, considering effective numerical methods becomes particularly necessary. Although the EM approximation \eqref{eq:5} provides a numerical approximation in the previous section, it will cost a lot of calculations on stochastic integrals. In order to reduce this cost, we now modify the approximation \eqref{eq:5} in this section.

\indent Under the same settings with the EM approximation \eqref{eq:5}, we can modify it by the left rectangle rule \cite{Kloeden1992, Liang2017} as
\begin{equation}\label{eq:10}
\begin{aligned}
Y(t) &= y_0 + \int_0^t F_0(t,\underline{s},\hat{Y}(s)) \mathrm{d}s + \int_0^t F_1(t,\underline{s},\hat{Y}(s)) \mathrm{d}s \\
&\quad + \int_0^t F_2(t,\underline{s},\hat{Y}(s)) \mathrm{d}W(s),
\end{aligned}
\end{equation}
where $\underline{s} = t_n$ for $s\in[t_n,t_{n+1})$ and the simple step process $\hat{Y}(t) = \sum_{n=0}^N Y(t_n)\cdot 1_{[t_n,t_{n+1})}(t)$. In the algorithm implementation, we adopt its discrete-time form
\begin{equation}\label{eq:11}
\begin{aligned}
Y_n := Y(t_n) &= y_0 + \sum_{j=0}^{n-1} F_0(t_n,t_j,Y_j)h + \sum_{j=0}^{n-1} F_1(t_n,t_j,Y_j)h \\
&\quad + \sum_{j=0}^{n-1} F_2(t_n,t_j,Y_j)\Delta W_j, \quad n = 1,2,\cdots, N
\end{aligned}
\end{equation}
with $Y_0 = y_0$, where $\Delta W_j = W(t_{j+1}) - W(t_j),~j = 0,1,\cdots,N-1$ denote the increments of Brownian motion $W(t)$. As an advantage, we only need to simulate these increments without computing other stochastic integrals, thus the calculations will be greatly reduced.

\subsection{Strong convergence}

\indent In order to explain the convergence of the EM scheme \eqref{eq:10}, we will first list some useful lemmas.

\begin{lemma}\label{le:4.1}
Let $a,b\in(0,1]$. Then for any $t\in[t_n,t_{n+1}]$, $n = 1,2,\cdots,N-1$, we have
\begin{align*}
& \int_0^{t_n} \big|(t - \underline{s})^{a-b} - (t_n - \underline{s})^{a-b}\big| \mathrm{d}s \leq C h^{1\wedge(1+a-b)}, \\
& \int_0^t \big|(t - \underline{s})^{a-b} - (t - s)^{a-b}\big| \mathrm{d}s \leq C h^{1\wedge(1+a-b)},
\end{align*}
where the positive constant $C$ depends on $T,a,b$, but not on $h$.
\end{lemma}

\begin{proof}
Setting $a,b\in(0,1]$ implies $a-b\in(-1,1)$. Firstly, for $a-b\in[0,1)$, one can derive
\begin{align*}
&\quad\ \int_0^{t_n} \big|(t - \underline{s})^{a-b} - (t_n - \underline{s})^{a-b}\big| \mathrm{d}s = \int_0^{t_n} (t - \underline{s})^{a-b} - (t_n - \underline{s})^{a-b} \mathrm{d}s \\
&\leq \sum_{i=0}^{n-1} \int_{t_i}^{t_{i+1}} (t_{n+1} -t_i)^{a-b} - (t_n -t_i)^{a-b} \mathrm{d}s \\
&= h \sum_{i=0}^{n-1} \left( \big( (n-i+1)h \big)^{a-b} - \big( (n-i)h \big)^{a-b} \right) \\
&= h \left( \big( (n+1)h \big)^{a-b} - h^{a-b} \right) \leq T^{a-b} h,
\end{align*}
and
\begin{align*}
&\quad \int_0^t \big|(t - \underline{s})^{a-b} - (t - s)^{a-b}\big| \mathrm{d}s = \int_0^t (t - \underline{s})^{a-b} - (t - s)^{a-b} \mathrm{d}s \\
&\leq \sum_{i=0}^{n-1} \int_{t_i}^{t_{i+1}} (t_{n+1} - t_i)^{a-b} - (t_n - t_{i+1})^{a-b} \mathrm{d}s + \int_{t_{n}}^t (t - \underline{s})^{a-b} - (t - s)^{a-b} \mathrm{d}s \\
&\leq h \sum_{i=0}^{n-1} \left( \big( (n-i+1)h \big)^{a-b} - \big( (n-i-1)h \big)^{a-b} \right) + \int_{t_{n}}^t (t - \underline{s})^{a-b} \mathrm{d}s \\
&= h \left( \big( (n+1)h \big)^{a-b} + (nh)^{a-b} - h^{a-b} \right) + \int_{t_{n}}^t (t - t_{n})^{a-b} \mathrm{d}s \\
&\leq 2T^{a-b} h + h^{1+a-b} \leq Ch.
\end{align*}

Secondly, for $a-b\in(-1,0)$, one can read
\begin{align*}
&\quad\ \int_0^{t_n} \big|(t - \underline{s})^{a-b} - (t_n - \underline{s})^{a-b}\big| \mathrm{d}s = \int_0^{t_n} (t_n - \underline{s})^{a-b} - (t - \underline{s})^{a-b} \mathrm{d}s \\
&\leq \sum_{i=0}^{n-1} \int_{t_i}^{t_{i+1}} (t_n -t_i)^{a-b} - (t_{n+1} -t_i)^{a-b} \mathrm{d}s \\
&= h^{1+a-b} \sum_{i=0}^{n-1} \left( (n-i)^{a-b} - (n-i+1)^{a-b} \right) \\
&= h^{1+a-b} \left( 1 - (n+1)^{a-b} \right) \leq h^{1+a-b},
\end{align*}
and
\begin{align*}
&\quad \int_0^t \big|(t - \underline{s})^{a-b} - (t - s)^{a-b}\big| \mathrm{d}s = \int_0^t (t - s)^{a-b} - (t - \underline{s})^{a-b} \mathrm{d}s \\
&\leq \sum_{i=0}^{n-2} \int_{t_i}^{t_{i+1}} (t_n - t_{i+1})^{a-b} - (t_{n+1} - t_i)^{a-b} \mathrm{d}s + \int_{t_{n-1}}^t (t - s)^{a-b} - (t - \underline{s})^{a-b} \mathrm{d}s \\
&\leq h^{1+a-b} \sum_{i=0}^{n-2} \left( (n-i-1)^{a-b} - (n-i+1)^{a-b} \right) + \int_{t_{n-1}}^t (t - s)^{a-b} \mathrm{d}s \\
&= h^{1+a-b} \left( 1 + 2^{a-b} - n^{a-b} - (n+1)^{a-b} \right) + \frac{(t-t_{n-1})^{(1+a-b)}}{1+a-b} \leq C h^{1+a-b}.
\end{align*}
The proof is now complete.
\end{proof}

\begin{lemma}\label{le:4.2}
Let $a\in(0,1]$, $b\in(0,\frac{1}{2})$. Then $a-b\in(-\frac{1}{2},1)$, and for any $t\in[t_n,t_{n+1}]$, $n = 1,2,\cdots,N-1$, we have
\begin{align*}
\int_0^{t_n} \big|(t - \underline{s})^{a-b} - (t_n - \underline{s})^{a-b}\big|^2 \mathrm{d}s \leq
\begin{cases}
C h^2, &\text{\quad if } a-b \in (\frac{1}{2},1), \\
C h^{2-\epsilon}, &\text{\quad if } a-b = \frac{1}{2}, \\
C h^{1+2(a-b)}, &\text{\quad if } a-b \in (-\frac{1}{2},\frac{1}{2}),
\end{cases}
\end{align*}
\begin{align*}
\int_0^t \big|(t - \underline{s})^{a-b} - (t - s)^{a-b}\big|^2 \mathrm{d}s \leq
\begin{cases}
C h^2, &\text{\quad if } a-b \in (\frac{1}{2},1), \\
C h^{2-\epsilon}, &\text{\quad if } a-b = \frac{1}{2}, \\
C h^{1+2(a-b)}, &\text{\quad if } a-b \in (-\frac{1}{2},\frac{1}{2}),
\end{cases}
\end{align*}
where the positive constant $C$ depends on $T,a,b,\epsilon$, but not on $h$, and $\epsilon \in (0,1)$.
\end{lemma}

\begin{proof}
Case 1: $a-b \in (\frac{1}{2},1)$. It follows from $a-b-1<0$ and $2(a-b-1)>-1$ that
\begin{align*}
&\quad\ \int_0^{t_n} \big| (t - \underline{s})^{a-b} - (t_n - \underline{s})^{a-b} \big|^2 \mathrm{d}s = \int_0^{t_n} \big| (a-b) \int_{t_n}^t (\tau - \underline{s})^{a-b-1} \mathrm{d}\tau \big|^2 \mathrm{d}s \\
&\leq (a-b)^2 \int_0^{t_n} \Big| \int_{t_n}^t (t_n - s)^{a-b-1} \mathrm{d}\tau \Big|^2 \mathrm{d}s \\
&\leq (a-b)^2 h^2 \int_0^{t_n} (t_n - s)^{2(a-b-1)} \mathrm{d}s \leq C h^2.
\end{align*}
Similarly,
\begin{align*}
\int_0^t \big| (t - \underline{s})^{a-b} - (t - s)^{a-b} \big|^2 \mathrm{d}s = \int_0^t \big| (a-b) \int_{\underline{s}}^s (t - \tau)^{a-b-1} \mathrm{d}\tau \big|^2 \mathrm{d}s \leq Ch^2.
\end{align*}

Case 2: $a-b = \frac{1}{2}$. For any $\epsilon \in (0,1)$, by Cauchy--Schwarz's inequality, one gets
\begin{align*}
&\quad\ \int_0^{t_n} \big| (t - \underline{s})^{a-b} - (t_n - \underline{s})^{a-b} \big|^2 \mathrm{d}s = \int_0^{t_n} \big| \frac{1}{2} \int_{t_n}^t (\tau - \underline{s})^{-\frac{1}{2}} \mathrm{d}\tau \big|^2 \mathrm{d}s \\
&\leq \frac{1}{4} \int_0^{t_n} \big| \int_{t_n}^t (\tau - s)^{-\frac{1}{2}} \mathrm{d}\tau \big|^2 \mathrm{d}s \leq \frac{1}{4} \int_0^{t_n} \int_{t_n}^t (\tau - s)^{-\epsilon} \mathrm{d}\tau \int_{t_n}^t (\tau - s)^{-1+\epsilon} \mathrm{d}\tau \mathrm{d}s \\
&= \frac{1}{4} \int_0^{t_n} \frac{1}{1-\epsilon} \big( (t-s)^{1-\epsilon} - (t_n-s)^{1-\epsilon} \big) \int_{t_n}^t (\tau - s)^{-1+\epsilon} \mathrm{d}\tau \mathrm{d}s \\
&\leq \frac{1}{4(1-\epsilon)} h^{1-\epsilon} \int_0^{t_n} \int_{t_n}^t (t_n - s)^{-1+\epsilon} \mathrm{d}\tau \mathrm{d}s \leq \frac{1}{4(1-\epsilon)} h^{2-\epsilon} \int_0^{t_n} (t_n - s)^{-1+\epsilon} \mathrm{d}s \\
&\leq \frac{T^{\epsilon}}{4(1-\epsilon)\epsilon} h^{2-\epsilon}.
\end{align*}
Analogously,
\begin{align*}
&\quad\ \int_0^t \big| (t - \underline{s})^{a-b} - (t - s)^{a-b} \big|^2 \mathrm{d}s = \int_0^{t} \big| \frac{1}{2} \int_{\underline{s}}^s (t - \tau)^{-\frac{1}{2}} \mathrm{d}\tau \big|^2 \mathrm{d}s \\
&\leq \frac{1}{4} \int_0^{t} \int_{\underline{s}}^s (t - \tau)^{-\epsilon} \mathrm{d}\tau \int_{\underline{s}}^s (t - \tau)^{-1+\epsilon} \mathrm{d}\tau \mathrm{d}s \\
&= \frac{1}{4} \int_0^{t} \frac{1}{1-\epsilon} \big( (t-\underline{s})^{1-\epsilon} - (t-s)^{1-\epsilon} \big) \int_{\underline{s}}^s (t - \tau)^{-1+\epsilon} \mathrm{d}\tau \mathrm{d}s \\
&\leq \frac{1}{4(1-\epsilon)} h^{1-\epsilon} \int_0^{t} \int_{\underline{s}}^s (t - s)^{-1+\epsilon} \mathrm{d}\tau \mathrm{d}s \leq \frac{1}{4(1-\epsilon)} h^{2-\epsilon} \int_0^{t} (t - s)^{-1+\epsilon} \mathrm{d}s \\
&\leq \frac{T^{\epsilon}}{4(1-\epsilon)\epsilon} h^{2-\epsilon}.
\end{align*}

Case 3: $a-b \in (-\frac{1}{2},\frac{1}{2})$. It follows from $a-b-1<0$ and $2(a-b-1)<-1$ that
\begin{align*}
&\quad\ \int_0^{t_n} \big|(t - \underline{s})^{a-b} - (t_n - \underline{s})^{a-b}\big|^2 \mathrm{d}s = \sum_{i=0}^{n-1} \int_{t_i}^{t_{i+1}} \big| (a-b) \int_{t_n}^t (\tau - t_i)^{a-b-1} \mathrm{d}\tau \big|^2 \mathrm{d}s \\
&\leq (a-b)^2 \sum_{i=0}^{n-1} \int_{t_i}^{t_{i+1}} \big| \int_{t_n}^t (t_n - t_i)^{a-b-1} \mathrm{d}\tau \big|^2 \mathrm{d}s \leq (a-b)^2 h^3 \sum_{i=0}^{n-1} \big( (n-i)h \big)^{2(a-b-1)} \\
&= (a-b)^2 h^{1+2(a-b)} \sum_{j=1}^n j^{2(a-b-1)} \leq C h^{1+2(a-b)}.
\end{align*}
Similarly,
\begin{align*}
&\quad\ \int_0^t \big|(t - \underline{s})^{a-b} - (t - s)^{a-b}\big|^2 \mathrm{d}s \\
&= \sum_{i=0}^{n-2} \int_{t_i}^{t_{i+1}} \big| (a-b) \int_{\underline{s}}^s (t - \tau)^{a-b-1} d\tau \big|^2 \mathrm{d}s + \int_{t_{n-1}}^t \big| (t - \underline{s})^{a-b} - (t - s)^{a-b} \big|^2 \mathrm{d}s \\
&\leq (a-b)^2 \sum_{i=0}^{n-2} \int_{t_i}^{t_{i+1}} \big| \int_{t_i}^s (t_n - t_{i+1})^{a-b-1} d\tau \big|^2 \mathrm{d}s \\
&\quad\ + \int_{t_{n-1}}^t 2 \big| (t - \underline{s})^{2(a-b)} + (t - s)^{2(a-b)} \big| \mathrm{d}s \\
&\leq (a-b)^2 h^3 \sum_{i=0}^{n-2} \big( (n-i-1)h \big)^{2(a-b-1)} + C h^{1+2(a-b)} \\
&\leq (a-b)^2 h^{1+2(a-b)} \sum_{j=1}^{n-1} j^{2(a-b-1)} + C h^{1+2(a-b)} \leq C h^{1+2(a-b)}.
\end{align*}
Hereto, Lemma \ref{le:4.2} has been proven.
\end{proof}

\begin{lemma}\label{le:4.3}
If the Assumption \ref{as:2.4} holds, then there exists a positive constant $C$ independent of $h$ such that for any integer $p \geq 2$,
\begin{align*}
\mathbb{E}\Big[|Y(t)|^p\Big] \leq C \quad \text{and} \quad \mathbb{E}\Big[|\hat{Y}(t)|^p\Big] \leq C, \qquad \forall\, t\in\mathcal{J}.
\end{align*}
\end{lemma}

\begin{proof}
The proof follows exactly the same lines as the proof of Lemma \ref{le:3.3}, so the details are omitted.
\end{proof}

\begin{lemma}\label{le:4.4}
If the Assumptions \ref{as:2.1} and \ref{as:2.4} hold, then there exists a positive constant $C$ independent of $h$ such that
\begin{align*}
\mathbb{E}\Big[|Y(t) - \hat{Y}(t)|^2\Big] \leq C h^{2\alpha}, \qquad \forall\, t\in\mathcal{J}.
\end{align*}
\end{lemma}

\begin{proof}
For arbitrary $t\in\mathcal{J}$, there is a unique integer $n$ such that $t\in[t_n,t_{n+1})$ and $\hat{Y}(t) = Y(t_n)$. Then, it follows from \eqref{eq:10} that
\begin{align*}
&\quad\ \mathbb{E}\Big[ |Y(t) - \hat{Y}(t)|^2\Big] = \mathbb{E}\Big[|Y(t) - Y(t_n)|^2 \Big] \\
&\leq 3 \bigg\{ \mathbb{E}\Big[ \big|\int_0^t F_0(t,\underline{s},\hat{Y}(s)) \mathrm{d}s - \int_0^{t_n} F_0(t_n,\underline{s},\hat{Y}(s)) \mathrm{d}s \big|^2 \Big] \\
&\qquad\ + \mathbb{E}\Big[ \big| \int_{0}^{t} F_1(t,\underline{s},\hat{Y}(s)) \mathrm{d}s - \int_{0}^{t_n} F_1(t_n,\underline{s},\hat{Y}(s)) \mathrm{d}s \big|^2 \Big] \\
&\qquad\ + \mathbb{E}\Big[ \big| \int_{0}^{t} F_2(t,\underline{s},\hat{Y}(s)) \mathrm{d}W(s) - \int_{0}^{t_n} F_2(t_n,\underline{s},\hat{Y}(s)) \mathrm{d}W(s) \big|^2 \Big] \bigg\} \\
&=: 3 \left\{ \widetilde{\mathcal{K}_0} + \widetilde{\mathcal{K}_1} + \widetilde{\mathcal{K}_2} \right\}.
\end{align*}
Applying H\"{o}lder's inequality and Assumptions \ref{as:2.1}, \ref{as:2.4} as well as Lemmas \ref{le:4.1} and \ref{le:4.3} to obtain
\begin{align*}
\widetilde{\mathcal{K}_0} &\leq 2 \bigg\{ \mathbb{E}\Big[ \big| \int_{0}^{t_n} \big( (t-\underline{s})^{\alpha-1} - (t_n-\underline{s})^{\alpha-1} \big) f_0(\underline{s},\hat{Y}(s)) \mathrm{d}s \big|^2 \Big] \\
&\qquad\ + \mathbb{E}\Big[ \big| \int_{t_n}^{t} (t-\underline{s})^{\alpha-1}f_0(\underline{s},\hat{Y}(s)) \mathrm{d}s \big|^2 \Big] \bigg\} \\
&\leq C \bigg\{ \int_0^{t_n} |(t-\underline{s})^{\alpha-1} - (t_n-\underline{s})^{\alpha-1}| \mathrm{d} s \cdot \\
&\qquad\qquad\qquad \int_0^{t_n} |(t-\underline{s})^{\alpha-1} - (t_n-\underline{s})^{\alpha-1}| \left(1+\mathbb{E}\Big[|\hat{Y}(s)|^2\Big]\right) \mathrm{d} s \\
&\qquad\ \ + \int_{t_n}^t (t-s)^{\alpha-1} \mathrm{d}s \cdot \int_{t_n}^t (t-s)^{\alpha-1} \left(1+\mathbb{E}\Big[|\hat{Y}(s)|^2\Big]\right) \mathrm{d} s \bigg\} \\
&\leq C h^{2\alpha}.
\end{align*}
Similarly, we can deal with $\widetilde{\mathcal{K}_1}$, and have
\begin{align*}
\widetilde{\mathcal{K}_1} \leq C h^{2 \wedge 2(1 + \alpha - \beta_1)} \leq C h^{2\alpha}.
\end{align*}
Now, we proceed to deal with $\widetilde{\mathcal{K}_2}$. Using H\"{o}lder's inequality, It\^{o} isometry and Assumptions \ref{as:2.1}, \ref{as:2.4} as well as Lemmas \ref{le:4.2} and \ref{le:4.3} to derive
\begin{align*}
\widetilde{\mathcal{K}_2} &\leq C \bigg\{ \mathbb{E}\Big[ \int_0^{t_n} \big( (t-\underline{s})^{\alpha-\beta_2} - (t_n-\underline{s})^{\alpha-\beta_2} \big)^2 \sup_{\underline{s}\leq u\leq t}\big|f_2(u,\underline{s},\hat{Y}(s))\big|^2 \mathrm{d}s \Big] \\
&\qquad\ \ + \mathbb{E}\Big[ \int_0^{t_n} (t_n-\underline{s})^{2(\alpha-\beta_2)} \sup_{0\leq u\leq 1} \big|f_2((t-\underline{s})u+\underline{s},\underline{s},\hat{Y}(s)) - \\
&\qquad\qquad\qquad\qquad\qquad\qquad\qquad\qquad\ \ f_2((t_n-\underline{s})u+\underline{s},\underline{s},\hat{Y}(s)) \big|^2 \mathrm{d}s \Big] \\
&\qquad\ \ + \mathbb{E}\Big[ \int_{t_n}^t (t-\underline{s})^{2(\alpha-\beta_2)} \sup_{\underline{s}\leq u\leq t}\big|f_2(u,\underline{s},\hat{Y}(s))\big|^2 \mathrm{d}s \bigg\} \\
&\leq C \bigg\{ \int_0^{t_n} \big( (t-\underline{s})^{\alpha-\beta_2} - (t_n-\underline{s})^{\alpha-\beta_2} \big)^2 \left(1+\mathbb{E}\Big[ |\hat{Y}(s)|^2 \Big]\right) \mathrm{d}s \\
&\qquad\ \ + h^2 \int_0^{t_n} (t_n-s)^{2(\alpha-\beta_2)} \left(1+\mathbb{E}\Big[ |\hat{Y}(s)|^2 \Big]\right) \mathrm{d}s \\
&\qquad\ \ + \int_{t_n}^t (t-s)^{2(\alpha-\beta_2)} \left(1+\mathbb{E}\Big[ |\hat{Y}(s)|^2 \Big]\right) \mathrm{d}s \bigg\} \\
&\leq
\begin{cases}
Ch^{2-\epsilon}, & \mbox{if } \alpha - \beta_2 = \frac{1}{2}, \\
Ch^{2 \wedge (1 + 2\alpha - 2\beta_2)}, & \mbox{otherwise}.
\end{cases}
\end{align*}
If $\alpha - \beta_2 = \frac{1}{2}$, then $2-2\alpha \in (0,1)$ since $\beta_2 \in (0,\frac{1}{2})$. Taking $\epsilon = 2-2\alpha$ yields
\begin{align*}
\widetilde{\mathcal{K}_2} \leq Ch^{2\alpha}.
\end{align*}
If $\alpha - \beta_2 \neq \frac{1}{2}$, then the assumption $\beta_2 \in (0,\frac{1}{2})$ implies
\begin{align*}
\widetilde{\mathcal{K}_2} \leq Ch^{2 \wedge (1 + 2\alpha - 2\beta_2)} \leq Ch^{2\alpha}.
\end{align*}
Hence,
\begin{align*}
\mathbb{E}\Big[|Y(t) - \hat{Y}(t)|^2\Big] \leq 3 \Big( \widetilde{\mathcal{K}_0} + \widetilde{\mathcal{K}_1} + \widetilde{\mathcal{K}_2} \Big) \leq C h^{2\alpha},
\end{align*}
which completes the proof.
\end{proof}

Now, we prove the mean-square convergence theorem of the EM method \eqref{eq:10}.

\begin{theorem}\label{th.4.5}
Under the Assumptions \ref{as:2.1}--\ref{as:2.4}, the EM solution $Y(t)$ defined by \eqref{eq:10} converges to the exact solution $y(t)$ in mean square sense, i.e.,
\begin{align*}
\lim_{h\rightarrow0}\mathbb{E}\Big[|Y(t) - y(t)|^2\Big] = 0, \qquad \forall\, t\in\mathcal{J}.
\end{align*}
\end{theorem}

\begin{proof}
For simplicity, let the error $e(t) = Y(t) - y(t)$ and for each integer $m\geq1$,
\begin{align*}
\tau_m = \inf \big\{t\geq0: |Y(t)|\geq m \big\}, \qquad \theta_m = \tau_m \wedge \rho_m,
\end{align*}
where $\rho_m$ has been given in \eqref{eq:7}. Similar to \eqref{eq:9}, it is follows from Lemma \ref{le:4.3} that for any $\delta > 0$,
\begin{align}\label{eq:12}
\mathbb{E}\Big[|e(t)|^2\Big] &\leq \mathbb{E}\Big[|e(t\wedge\theta_m)|^2\Big] + \frac{2^{p+1}\delta M}{p} + \frac{(p-2)2M}{p\delta^{\frac{2}{p-2}}m^p},
\end{align}
here and in the rest of the proof, $M$ represents a positive constant independent of $m$, $\delta$ and $h$. Next, we focus on estimating the first term on the right side of \eqref{eq:12}. According to H\"{o}lder's inequality, it follows from the formulae \eqref{eq:4} and \eqref{eq:10} that
\begin{align}\label{eq:13}
\mathbb{E}\Big[|e(t\wedge\theta_m)|^2\Big] &\leq 6\mathbb{E}\Bigg\{\Big| \int_0^{t\wedge\theta_m} F_0(t\wedge\theta_m,\underline{s},\hat{Y}(s)) - F_0(t\wedge\theta_m,s,\hat{Y}(s)) \mathrm{d}s \Big|^2 \nonumber\\
&\qquad\ \ + \Big| \int_0^{t\wedge\theta_m} F_0(t\wedge\theta_m,s,\hat{Y}(s)) - F_0(t\wedge\theta_m,s,y(s)) \mathrm{d}s \Big|^2 \nonumber\\
&\qquad\ \ + \Big| \int_0^{t\wedge\theta_m} F_1(t\wedge\theta_m,\underline{s},\hat{Y}(s)) - F_1(t\wedge\theta_m,s,\hat{Y}(s)) \mathrm{d}s \Big|^2 \nonumber\\
&\qquad\ \ + \Big| \int_0^{t\wedge\theta_m} F_1(t\wedge\theta_m,s,\hat{Y}(s)) - F_1(t\wedge\theta_m,s,y(s)) \mathrm{d}s \Big|^2 \nonumber\\
&\qquad\ \ + \Big| \int_0^{t\wedge\theta_m} F_2(t\wedge\theta_m,\underline{s},\hat{Y}(s)) - F_2(t\wedge\theta_m,s,\hat{Y}(s)) \mathrm{d}W(s) \Big|^2 \nonumber\\
&\qquad\ \ + \Big| \int_0^{t\wedge\theta_m} F_2(t\wedge\theta_m,s,\hat{Y}(s)) - F_2(t\wedge\theta_m,s,y(s)) \mathrm{d}W(s) \Big|^2 \Bigg\} \nonumber\\
&=: 6 \Big\{ \mathcal{H}_1 + \mathcal{H}_2 + \mathcal{H}_3 + \mathcal{H}_4 + \mathcal{H}_5 + \mathcal{H}_6 \Big\}.
\end{align}
Applying Cauchy--Schwarz's inequality and It\^{o} isometry as well as Assumptions \ref{as:2.1}, \ref{as:2.2}, \ref{as:2.4} and Lemmas \ref{le:4.1}--\ref{le:4.3}, it follows from the derivation steps of $\widetilde{\mathcal{K}_0}$, $\widetilde{\mathcal{K}_1}$ and $\widetilde{\mathcal{K}_2}$ (see the proof of Lemma \ref{le:4.4}) that
\begin{align}\label{eq:14}
\mathcal{H}_1 + \mathcal{H}_3 + \mathcal{H}_5 \leq C h^{2\alpha}.
\end{align}
By Cauchy--Schwarz's inequality, Assumption \ref{as:2.3} and It\^{o} isometry, it follows from the derivation steps of \eqref{eq:6} that
\begin{align}\label{eq:15}
\mathcal{H}_2 + \mathcal{H}_4 + \mathcal{H}_6 \leq CK_m^2 \int_0^{t\wedge\theta_m} (t\wedge\theta_m-s)^{\alpha-1} \mathbb{E}\Big[ |\hat{Y}(s) - Y(s)|^2 + |e(s\wedge\theta_m)|^2 \Big] \mathrm{d}s.
\end{align}
Now, combining \eqref{eq:13}--\eqref{eq:15}, and applying Lemma \ref{le:4.4} as well as weakly singular Gronwall's inequality \cite[Corollary 2]{YeGao2007}, we arrive at
\begin{align}\label{eq:16}
\mathbb{E}\Big[|e(t\wedge\theta_m)|^2\Big] \leq C_m h^{2\alpha},
\end{align}
where the positive constant $C_m$ depends on $m$, but not on $h$ and $\delta$. Inserting \eqref{eq:16} into \eqref{eq:12} yields
\begin{align*}
\mathbb{E}\Big[|e(t)|^2\Big] \leq C_m h^{2\alpha} + \frac{2^{p+1}\delta M}{p} + \frac{(p-2)2M}{p\delta^{\frac{2}{p-2}}m^p}.
\end{align*}
Therefore, for any given $\varepsilon > 0$, one can choose $\delta$ and $m$ such that
\begin{align*}
\frac{2^{p+1}\delta M}{p} < \frac{\varepsilon}{3}, \qquad \frac{(p-2)2M}{p\delta^{\frac{2}{p-2}}m^p} < \frac{\varepsilon}{3},
\end{align*}
and then $h$ can be taken sufficiently small such that $C_m h^{2\alpha} < \frac{\varepsilon}{3}$. As a result,
\begin{align*}
\lim_{h\rightarrow0}\mathbb{E}\Big[|e(t)|^2\Big] = 0, \qquad \forall\, t\in\mathcal{J}.
\end{align*}
The proof is complete.
\end{proof}

\begin{remark}
It is a classical result \cite{Jentzen2011divergen} that the EM method is divergent for SDEs with superlinearly growing coefficients in mean square sense, which motivates the use of the linear growth assumption in the above theorem.
\end{remark}

\subsection{Convergence rate}

\indent The convergence rate of the numerical scheme can display its computational efficiency, so we give the following theorem.

\begin{theorem}\label{th.4.6}
If the assumptions of Theorem \ref{th.4.5} and the global Lipschitz condition \eqref{eq:2} are satisfied, then there exists a positive constant $C$ independent of $h$ such that
\begin{align*}
\bigg( \mathbb{E}\Big[|Y(t) - y(t)|^2\Big] \bigg)^{\frac{1}{2}} \leq C h^{\alpha}, \qquad \forall\, t\in\mathcal{J}.
\end{align*}
\end{theorem}

\begin{proof}
It follows from \eqref{eq:4} and \eqref{eq:10} that
\begin{align}\label{eq:17}
\mathbb{E}\Big[|Y(t) - y(t)|^2\Big] &\leq 6 \bigg\{ \mathbb{E}\Big[ \big| \int_0^t F_0(t,\underline{s},\hat{Y}(s)) - F_0(t,s,\hat{Y}(s)) \mathrm{d}s \big|^2 \Big] \nonumber\\
&\qquad\ + \mathbb{E}\Big[ \big| \int_0^t F_0(t,s,\hat{Y}(s)) - F_0(t,s,y(s)) \mathrm{d}s \big|^2 \Big] \nonumber\\
&\qquad\ + \mathbb{E}\Big[ \big| \int_0^t F_1(t,\underline{s},\hat{Y}(s)) - F_1(t,s,\hat{Y}(s)) \mathrm{d}s \big|^2 \Big] \nonumber\\
&\qquad\ + \mathbb{E}\Big[ \big| \int_0^t F_1(t,s,\hat{Y}(s)) - F_1(t,s,y(s)) \mathrm{d}s \big|^2 \Big] \nonumber\\
&\qquad\ + \mathbb{E}\Big[ \big| \int_0^t F_2(t,\underline{s},\hat{Y}(s)) - F_2(t,s,\hat{Y}(s)) \mathrm{d}W(s) \big|^2 \Big] \nonumber\\
&\qquad\ + \mathbb{E}\Big[ \big| \int_0^t F_2(t,s,\hat{Y}(s)) - F_2(t,s,y(s)) \mathrm{d}W(s) \big|^2 \Big] \bigg\} \nonumber\\
&=: 6 \Big\{ \widetilde{\mathcal{H}_1} + \widetilde{\mathcal{H}_2} + \widetilde{\mathcal{H}_3} + \widetilde{\mathcal{H}_4} + \widetilde{\mathcal{H}_5} + \widetilde{\mathcal{H}_6} \Big\}.
\end{align}
Using Cauchy--Schwarz's inequality and It\^{o} isometry as well as Assumptions \ref{as:2.1}, \ref{as:2.2}, \ref{as:2.4} and Lemmas \ref{le:4.1}--\ref{le:4.3}, it follows from similar derivation steps of $\widetilde{\mathcal{K}_0}$, $\widetilde{\mathcal{K}_1}$ and $\widetilde{\mathcal{K}_2}$ (see the proof of Lemma \ref{le:4.4}) that
\begin{align}\label{eq:18}
\widetilde{\mathcal{H}_1} + \widetilde{\mathcal{H}_3} + \widetilde{\mathcal{H}_5} \leq C h^{2\alpha}.
\end{align}
By Cauchy--Schwarz's inequality, global Lipschitz condition \eqref{eq:2} and It\^{o} isometry, it follows from the derivation steps of \eqref{eq:6} and the triangle inequality that
\begin{align}\label{eq:19}
\widetilde{\mathcal{H}_2} + \widetilde{\mathcal{H}_4} + \widetilde{\mathcal{H}_6} \leq CK^2 \int_0^{t} (t-s)^{\alpha-1} \mathbb{E}\Big[ |\hat{Y}(s) - Y(s)|^2 + |Y(t) - y(t)|^2 \Big] \mathrm{d}s.
\end{align}
Finally, combining \eqref{eq:17}--\eqref{eq:19}, and applying Lemma \ref{le:4.4} as well as weakly singular Gronwall's inequality \cite[Corollary 2]{YeGao2007} shows that
\begin{align*}
\mathbb{E}\Big[|Y(t) - y(t)|^2\Big] \leq C h^{2\alpha},
\end{align*}
which completes the proof.
\end{proof}

\begin{remark}
For the obtained convergence rate, we consider the following two special cases.
\begin{itemize}
  \item Case 1: When $\alpha = 1$, the SFIDE \eqref{eq:1} becomes the integer-order stochastic integro-differential equation, the EM method \eqref{eq:10} can attain strong first-order superconvergence, which actually improves the corresponding result of \cite[the case $H \in (0,\frac{1}{2})$ of Theorem 3.9]{YangYangYao2021}.
  \item Case 2: When $f_1=f_2\equiv0$, the SFIDE \eqref{eq:1} degenerates the deterministic fractional differential equation, and the obtained convergence rate is consistent with the existing results, e.g., \cite[Page 7]{IzzoMessina2018}.
\end{itemize}
Therefore, one can conclude that the result obtained in Theorem \ref{th.4.6} is very sharp.
\end{remark}

\begin{remark}
Under the global Lipschitz condition and linear growth condition, Theorem \ref{th.4.6} attains mean-square convergence rate of the explicit EM method \eqref{eq:10} for SFIDE \eqref{eq:1}. Recalling numerical methods of SDEs with superlinearly growing coefficients, we know that the explicit EM method fails to converge strongly to its exact solution \cite{Jentzen2011divergen}. In this case, implicit methods \cite{Higham2002} or explicitly tamed methods \cite{Jentzentamed2012} are still effective, the later of which is less computationally expensive. As a result, the tamed EM method may be a well candidate numerical method for SFIDE \eqref{eq:1} with superlinearly growing coefficients, which will be our next goal in the future.
\end{remark}

\section{Numerical experiments}
\label{sec.5}

\indent In this section, we verify the convergence rate of the EM method given in \eqref{eq:11} for some SFIDEs with weakly singular kernels. In a similar way as \cite{Cao2015}, we use sample average to approximate the expectation. More precisely, we measure the mean square error of numerical solutions by
\begin{align*}
\epsilon_{h,T} = \left( \frac{1}{5000}\sum_{i=1}^{5000}\left| Y_{h}(T,\omega_i) - Y_{\frac{h}{2}}(T,\omega_i) \right|^2 \right)^{1/2},
\end{align*}
where $\omega_i$ denotes the $i$th single sample path.

\begin{example}\label{ex:5.1}
Consider the $1$-dimensional SFIDE with $r = 1$
\begin{align*}
D^{\alpha}y(t) &= \sin(ty(t)) + \int_0^t \frac{ts\cos(y(s))}{(t-s)^{\beta_1}} \mathrm{d}s + \int_0^t \frac{ts\cos(y(s))}{(t-s)^{\beta_2}} \mathrm{d}W(s)
\end{align*}
for $t\in[0,1]$ and the initial value $y(0) = 1$.

\begin{figure}[htp]
\begin{center}
  \includegraphics[width=5in]{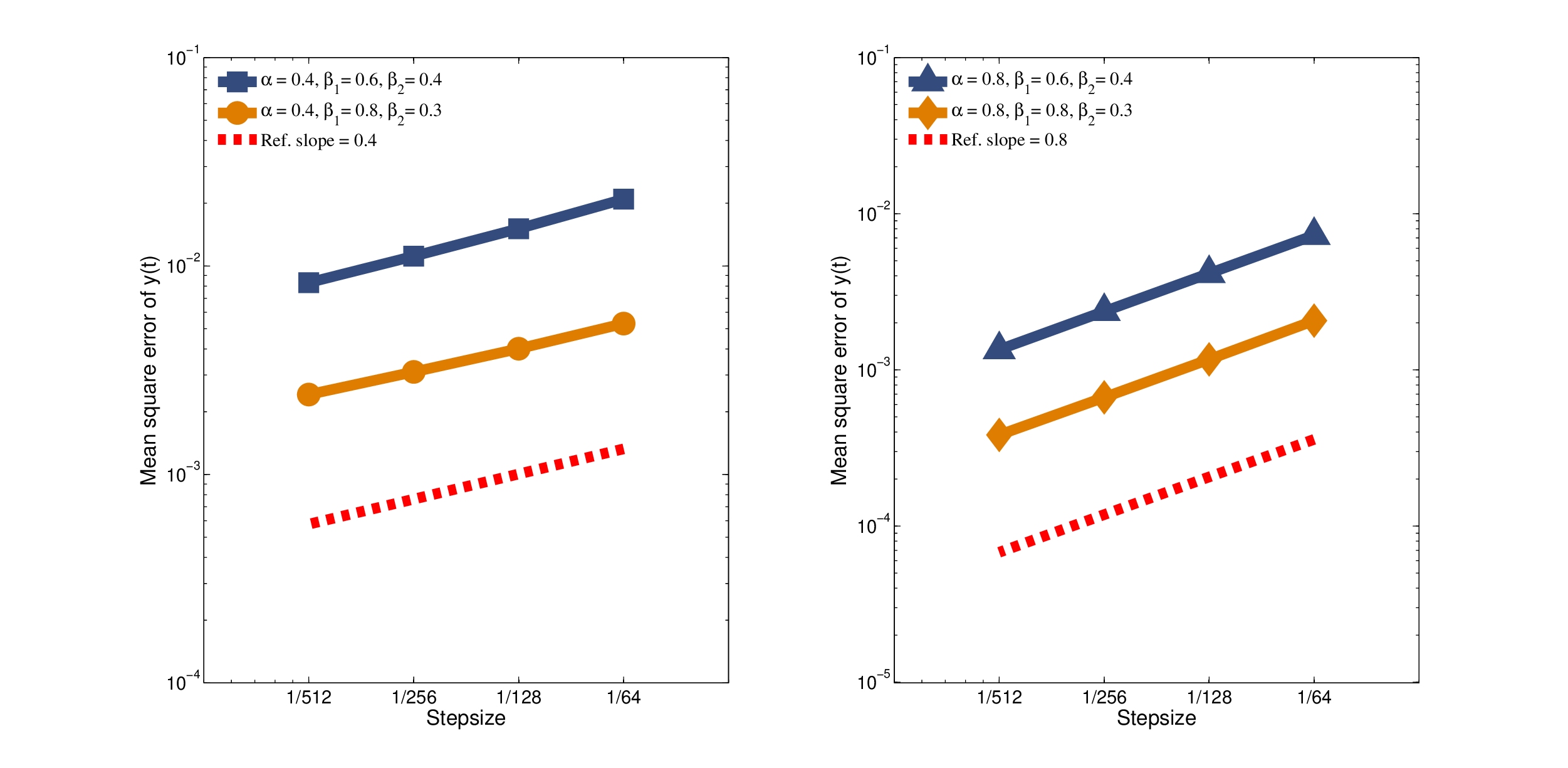}
  \caption{The mean square errors of the EM scheme \eqref{eq:11} for Example \ref{ex:5.1}.}\label{fig:1}
\end{center}
\end{figure}

\indent We can verify that the functions $f_i~(i=0,1,2)$ satisfy the hypotheses of Theorem \ref{th.4.6}. The computing results are shown in Figure \ref{fig:1}. As shown in Figure \ref{fig:1}, the convergence rate is $\alpha$, and the arguments $\alpha$, $\beta_1$ and $\beta_2$ will affect the error constant.
\end{example}

\begin{example}\label{ex:5.2}
Consider the $2$-dimensional SFIDE with $r = 2$
\begin{equation*}
\begin{array}{llll}
&D^{\alpha}\begin{bmatrix}
y_1(t) \\
y_2(t)
\end{bmatrix}
= \begin{bmatrix}
\sin(ty_2(t)) \\
ty_1(t)
\end{bmatrix}
+\begin{bmatrix}
\int_0^t \frac{s\sin(y_1(s)+y_2(s))}{(t-s)^{\beta_1}} \mathrm{d}s + \int_0^t \frac{sy_1(s)}{(t-s)^{\beta_2}} \mathrm{d}W_1(s) \\
\int_0^t \frac{s\cos(y_1(s)+y_2(s))}{(t-s)^{\beta_1}} \mathrm{d}s + \int_0^t \frac{s\cos(y_2(s))}{(t-s)^{\beta_2}} \mathrm{d}W_1(s)
\end{bmatrix} \\
&\qquad\qquad\qquad\ +\begin{bmatrix}
\int_0^t \frac{s\cos(y_1(s)+y_2(s))}{(t-s)^{\beta_2}} \mathrm{d}W_2(s) \\
\int_0^t \frac{s\sin(y_1(s)+y_2(s))}{(t-s)^{\beta_2}} \mathrm{d}W_2(s)
\end{bmatrix}
\end{array}
\end{equation*}
for $t\in[0,1]$ and the initial value $(y_1(0),y_2(0))^T = (0,0)^T$.

\begin{figure}[htp]
\begin{center}
  \includegraphics[width=5in]{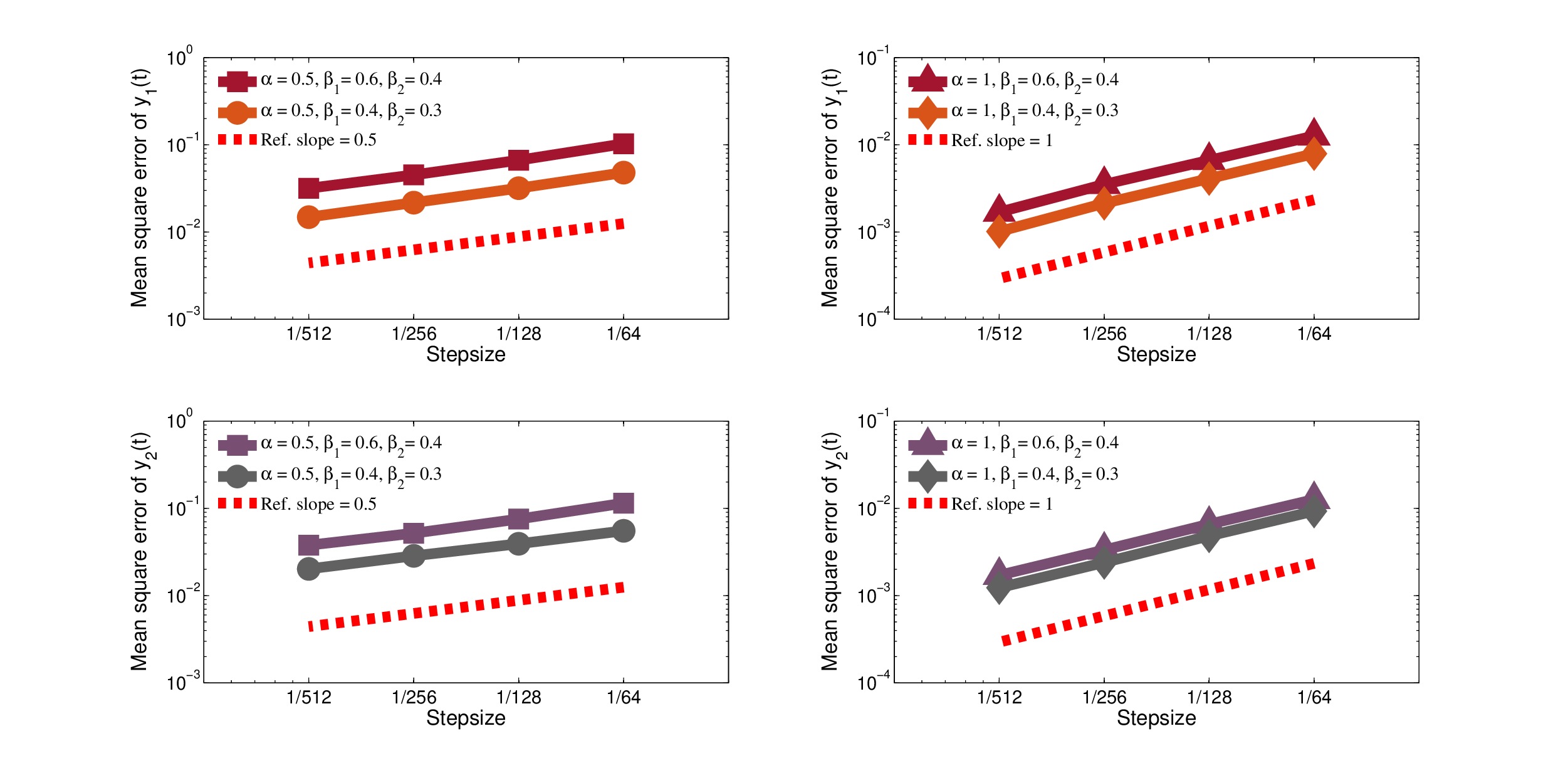}
  \caption{The mean square errors of the EM scheme \eqref{eq:11} for Example \ref{ex:5.2}.}\label{fig:2}
\end{center}
\end{figure}

\indent The functions $f_i~(i=0,1,2)$ can be verified to satisfy the hypotheses of Theorem \ref{th.4.6}. The corresponding calculation results are displayed in Figure\ \ref{fig:2}. As shown in Figure \ref{fig:2}, the same phenomenon appears as in Example \ref{ex:5.1}.
\end{example}

\section{Conclusion}
\label{sec.6}

\indent Under some relaxed conditions (e.g., the local Lipschitz condition), we obtain the existence, uniqueness and continuous dependence on the initial value (in mean square sense) of the true solution to SFIDEs with weakly singular kernels. Moreover, we analyze the strong convergence of the EM method \eqref{eq:10}, whose computational cost on stochastic integrals can be reduced. All of the numerical results are in line with our theoretical results.

\section*{Acknowledgments} Deep thanks go to the editor for his help and insightful comments, and the referees for their very careful reading of the manuscript and their helpful comments and suggestions, which greatly improved the quality of this article.


\medskip
Received xxxx 20xx; revised xxxx 20xx.
\medskip


\begin{thebibliography}{99}

\bibitem{Aghajani2012} (MR2872111) [10.2478/s13540-012-0005-4]
     \newblock A. Aghajani, Y. Jalilian, J. Trujillo,
     \newblock On the existence of solutions of fractional integro-differential equations,
     \newblock \emph{Fract. Calc. Appl. Anal.}, \textbf{15} (2012), 44--69.

\bibitem{AnhDoanHuong2019} (MR3873977) [10.1016/j.spl.2018.10.010]
    \newblock P.T. Anh, T.S. Doan, P.T. Huong,
    \newblock A variation of constant formula for Caputo fractional stochastic differential equations,
    \newblock \emph{Statist. Probab. Lett.}, \textbf{145} (2019), 351--358.

\bibitem{Asgari2014} (MR3263911) [10.5899/2014/cna-00212]
     \newblock M. Asgari,
     \newblock Block pulse approximation of fractional stochastic integro-differential equation,
     \newblock \emph{Commun. Numer. Anal.}, \textbf{2014} (2014), 1--7.

\bibitem{Badr2012} (MR2912586) [10.1155/2012/709106 ]
     \newblock A.A. Badr, H.S. El-Hoety,
     \newblock Monte--Carlo Galerkin approximation of fractional stochastic integro-differential equation,
     \newblock \emph{Math. Probl. Eng.}, \textbf{2012} (2012), 709106.

\bibitem{Balasubramaniam2017} (MR3673709) [10.1007/s10957-016-0865-6]
     \newblock P. Balasubramaniam, P. Tamilalagan,
     \newblock The solvability and optimal controls for impulsive fractional stochastic integro-differential equations via resolvent operators,
     \newblock \emph{J. Optim. Theory Appl.}, \textbf{174} (2017), 139--155.

\bibitem{Cao2015} (MR3305372) [10.1137/130942024]
     \newblock W. Cao, Z. Zhang, G.E. Karniadakis,
     \newblock Numerical methods for stochastic delay differential equations via the Wong--Zakai approximation,
     \newblock \emph{SIAM J. Sci. Comput.}, \textbf{37} (2015), A295--A318.

\bibitem{ChenKim2015} (MR3310354) [10.1016/j.spa.2014.11.005]
     \newblock Z.-Q. Chen, K.-H. Kim, P. Kim,
     \newblock Fractional time stochastic partial differential equations,
     \newblock \emph{Stochastic Process. Appl.}, \textbf{125} (2015), 1470--1499.

\bibitem{ContTankov2004} (MR2042661) [10.1201/9780203485217]
     \newblock R. Cont, P. Tankov,
     \newblock \emph{Financial Modelling with Jump Processes},
     \newblock Chapman and Hall/CRC, 2004.

\bibitem{DaPrato2014} (MR3236753) [10.1017/CBO9781107295513]
     \newblock G. Da Prato, J. Zabczyk,
     \newblock \emph{Stochastic Equations in Infinite Dimensions},
     \newblock Cambridge University Press, 2014.

\bibitem{Dai2019} (MR3921148) [10.1016/j.cam.2019.02.002]
     \newblock X. Dai, W. Bu, A. Xiao,
     \newblock Well-posedness and EM approximations for non-Lipschitz stochastic fractional integro-differential equations,
     \newblock \emph{J. Comput. Appl. Math.}, \textbf{356} (2019), 377--390.

\bibitem{Diethelm2010} (MR2680847) []
    \newblock K. Diethelm,
    \newblock \emph{The Analysis of Fractional Differential Equations:\ An Application-Oriented Exposition Using Differential Operators of Caputo Type},
    \newblock Springer, 2010.

\bibitem{DoanHuongKloedenVu2020} (MR4105671) [10.1016/j.cam.2020.112989]
    \newblock T.S. Doan, P.T. Huong, P.E. Kloeden, A.M. Vu,
    \newblock Euler--Maruyama scheme for Caputo stochastic fractional differential equations,
    \newblock \emph{J. Comput. Appl. Math.}, \textbf{380} (2020), 112989.

\bibitem{Higham2002} (MR1949404) [10.1137/s0036142901389530]
     \newblock D.J. Higham, X. Mao, A.M. Stuart,
     \newblock Strong convergence of Euler-type methods for nonlinear stochastic differential equations,
     \newblock \emph{SIAM J. Numer. Anal.}, \textbf{40} (2002), 1041--1063.

\bibitem{Jentzen2011divergen} (MR2795791) [10.1098/rspa.2010.0348]
     \newblock M. Hutzenthaler, A. Jentzen, P.E. Kloeden,
     \newblock Strong and weak divergence in finite time of Euler's method for stochastic differential equations with non-globally Lipschitz continuous coefficients,
     \newblock \emph{Proc. R. Soc. Lond. Ser. A Math. Phys. Eng. Sci.}, \textbf{467} (2011), 1563--1576.

\bibitem{Jentzentamed2012} (MR2985171) [10.1214/11-AAP803]
    \newblock M. Hutzenthaler, A. Jentzen, P.E. Kloeden,
    \newblock Strong convergence of an explicit numerical method for SDEs with nonglobally Lipschitz continuous coefficients,
    \newblock \emph{Ann. Appl. Probab.} \textbf{22} (2012), 1611--1641.

\bibitem{IzzoMessina2018} (MR3800780) [10.1007/s00009-018-1149-1]
    \newblock G. Izzo, E. Messina, A. Vecchio,
    \newblock Stability of numerical solutions for Abel--Volterra integral equations of the second kind,
    \newblock \emph{Mediterr. J. Math.}, \textbf{15} (2018), 113.

\bibitem{JinYan2019} (MR3978476) [10.1051/m2an/2019025]
     \newblock B. Jin, Y. Yan, Z. Zhou,
     \newblock Numerical approximation of stochastic time-fractional diffusion,
     \newblock \emph{ESAIM Math. Model. Numer. Anal.}, \textbf{53} (2019), 1245--1268.

\bibitem{Kamrani2015} (MR3296700) [10.1007/s11075-014-9839-7]
     \newblock M. Kamrani,
     \newblock Numerical solution of stochastic fractional differential equations,
     \newblock \emph{Numer. Algorithms}, \textbf{68} (2015), 81--93.

\bibitem{Kamrani2016} () [10.1016/j.ijleo.2016.07.087]
     \newblock M. Kamrani,
     \newblock Convergence of Galerkin method for the solution of stochastic fractional integro differential equations,
     \newblock \emph{Optik Int. J. Light Electron Opt.}, \textbf{127} (2016), 10049--10057.

\bibitem{Kloeden1992} (MR1214374) [10.1007/978-3-662-12616-5]
     \newblock P.E. Kloeden, E. Platen,
     \newblock \emph{Numerical Solution of Stochastic Differential Equations},
     \newblock Springer, 1992.

\bibitem{Lakshmikantham1995} (MR1336142) []
     \newblock V. Lakshmikantham, M. Rama Mohana Rao,
     \newblock \emph{Theory of Integro-Differential Equations},
     \newblock CRC Press, 1995.

\bibitem{Levin1960} () [10.1512/iumj.1960.9.59020]
     \newblock J.J. Levin, J.A. Nohel,
     \newblock On a system of integrodifferential equations occuring in reactor dynamics,
     \newblock \emph{J. Math. Mech.}, \textbf{9} (1960), 347--368.

\bibitem{LiLiu2017} (MR3704863) [10.1007/s10955-017-1866-z]
     \newblock L. Li, J.-G. Liu, J. Lu,
     \newblock Fractional stochastic differential equations satisfying fluctuation-dissipation theorem,
     \newblock \emph{J. Stat. Phys.}, \textbf{169} (2017), 316--339.

\bibitem{LiWang2019} (MR3912691) [10.1016/j.jde.2018.09.009]
     \newblock Y. Li, Y. Wang,
     \newblock The existence and asymptotic behavior of solutions to fractional stochastic evolution equations with infinite delay,
     \newblock \emph{J. Differential Equations}, \textbf{266} (2019), 3514--3558.

\bibitem{Liang2017} (MR3606089) [10.1016/j.cam.2016.11.005]
     \newblock H. Liang, Z. Yang, J. Gao,
     \newblock Strong superconvergence of the Euler--Maruyama method for linear stochastic Volterra integral equations,
     \newblock \emph{J. Comput. Appl. Math.}, \textbf{317} (2017), 447--457.

\bibitem{Maleki2015} (MR3353067) [10.1016/j.apm.2014.12.045]
     \newblock M. Maleki, M.T. Kajani,
     \newblock Numerical approximations for Volterra's population growth model with fractional order via a multi-domain pseudospectral method,
     \newblock \emph{Appl. Math. Model.}, \textbf{39} (2015), 4300--4308.

\bibitem{Mao2008} (MR2380366) [10.1533/9780857099402]
     \newblock X. Mao,
     \newblock \emph{Stochastic Differential Equations and Applications},
     \newblock Elsevier, 2008.
	
\bibitem{McKinley2018} (MR3856208) [10.1137/17m115517x]
     \newblock S.A. McKinley, H.D. Nguyen,
     \newblock Anomalous diffusion and the generalized Langevin equation,
     \newblock \emph{SIAM J. Math. Anal.}, \textbf{50} (2018), 5119--5160.

\bibitem{Mirzaee2017} () [10.1016/j.ijleo.2016.12.029]
     \newblock F. Mirzaee, N. Samadyar,
     \newblock Application of orthonormal Bernstein polynomials to construct a efficient scheme for solving fractional stochastic integro-differential equation,
     \newblock \emph{Optik Int. J. Light Electron Opt.}, \textbf{132} (2017), 262--273.

\bibitem{Mirzaee2018} (MR3913709) [10.1016/j.enganabound.2018.05.006]
     \newblock F. Mirzaee, N. Samadyar,
     \newblock On the numerical solution of fractional stochastic integro-differential equations via meshless discrete collocation method based on radial basis functions,
     \newblock \emph{Eng. Anal. Bound. Elem.}, \textbf{100} (2019), 246--255.
	
\bibitem{Mohammadi2015} (MR3573597) [10.5269/bspm.v35i1.28262]
     \newblock F. Mohammadi,
     \newblock Efficient Galerkin solution of stochastic fractional differential equations using second kind Chebyshev wavelets,
     \newblock \emph{Bol. Soc. Parana. Mat.}, \textbf{35} (2015), 195--215.

\bibitem{Momani2000} (MR1793439)
     \newblock S.M. Momani,
     \newblock Local and global existence theorems on fractional integro-differential equations,
     \newblock \emph{J. Fract. Calc.}, \textbf{18} (2000), 81--86.

\bibitem{Pedjeu2012} (MR2881656) [10.1016/j.chaos.2011.12.009]
     \newblock J.C. Pedjeu, G.S. Ladde,
     \newblock Stochastic fractional differential equations:\ Modeling, method and analysis,
     \newblock \emph{Chaos Solitons Fractals}, \textbf{45} (2012), 279--293.

\bibitem{Rao1975} (MR0358994) [10.1016/s0019-9958(75)90074-1]
     \newblock A.N.V. Rao, C.P. Tsokos,
     \newblock On the existence, uniqueness, and stability behavior of a random solution to a nonlinear perturbed stochastic integro-differential equation,
     \newblock \emph{Information and Control}, \textbf{27} (1975), 61--74.

\bibitem{Scudo1971} (MR0408866) [10.1016/0040-5809(71)90002-5]
    \newblock F.M. Scudo,
    \newblock Vito Volterra and theoretical ecology,
    \newblock \emph{Theoret. Population Biol.}, \textbf{2} (1971), 1--23.

\bibitem{SonHuong2018} (MR3854535) [10.1080/07362994.2018.1440243]
     \newblock D.T. Son, P.T. Huong, P.E. Kloeden, H.T. Tuan,
     \newblock Asymptotic separation between solutions of Caputo fractional stochastic differential equations,
     \newblock \emph{Stoch. Anal. Appl.}, \textbf{36} (2018), 654--664.

\bibitem{Taheri2017} (MR3634940) [10.1016/j.cam.2017.02.027]
     \newblock Z. Taheri, S. Javadi, E. Babolian,
     \newblock Numerical solution of stochastic fractional integro-differential equation by the spectral collocation method,
     \newblock \emph{J. Comput. Appl. Math.}, \textbf{321} (2017), 336--347.

\bibitem{Tarasov2009} () [10.1007/s11232-009-0029-z]
     \newblock V.E. Tarasov,
     \newblock Fractional integro-differential equations for electromagnetic waves in dielectric media,
     \newblock \emph{Theoret. and Math. Phys.}, \textbf{158} (2009), 355--359.

\bibitem{TeBeest1997} (MR1469945) [10.1137/s0036144595294850]
     \newblock K.G. TeBeest,
     \newblock Classroom Note:\ Numerical and analytical solutions of Volterra's population model,
     \newblock \emph{SIAM Rev.}, \textbf{39} (1997), 484--493.

\bibitem{Tien2013} (MR2968995) [10.1016/j.jmaa.2012.07.062]
     \newblock D.N. Tien,
     \newblock Fractional stochastic differential equations with applications to finance,
     \newblock \emph{J. Math. Anal. Appl.}, \textbf{397} (2013), 334--348.

\bibitem{Tuan2020} (MR4203018) [10.3934/dcdsb.2020318]
    \newblock H.T. Tuan,
    \newblock On the asymptotic behavior of solutions to time-fractional elliptic equations driven by a multiplicative white noise,
    \newblock \emph{Discrete Contin. Dyn. Syst. Ser. B}, \textbf{26} (2021), 1749--1762.

\bibitem{Wang2008} (MR2422961) [10.1016/j.spl.2007.10.007]
     \newblock Z. Wang,
     \newblock Existence and uniqueness of solutions to stochastic Volterra equations with singular kernels and non-Lipschitz coefficients,
     \newblock \emph{Statist. Probab. Lett.}, \textbf{78} (2008), 1062--1071.

\bibitem{WangXu2016} (MR3473117) [10.1016/j.na.2016.01.020]
     \newblock Y. Wang, J. Xu, P.E. Kloeden,
     \newblock Asymptotic behavior of stochastic lattice systems with a Caputo fractional time derivative,
     \newblock \emph{Nonlinear Anal.}, \textbf{135} (2016), 205--222.

\bibitem{YangYangYao2021} (MR4140608) [10.1016/j.cam.2020.113156]
    \newblock Z. Yang, H. Yang, Z. Yao,
    \newblock Strong convergence analysis for Volterra integro-differential equations with fractional Brownian motions,
    \newblock \emph{J. Comput. Appl. Math.}, \textbf{383} (2021), 113156.

\bibitem{YeGao2007} (MR2290034) [10.1016/j.jmaa.2006.05.061]
     \newblock H. Ye, J. Gao, Y. Ding,
     \newblock A generalized Gronwall inequality and its application to a fractional differential equation,
     \newblock \emph{J. Math. Anal. Appl.}, \textbf{328} (2007), 1075--1081.

\bibitem{Yuzbasi2013} (MR3020487) [10.1016/j.apm.2012.07.041]
     \newblock \c{S}. Y\"{u}zba\c{s}{\i},
     \newblock A numerical approximation for Volterra's population growth model with fractional order,
     \newblock \emph{Appl. Math. Model.}, \textbf{37} (2013), 3216--3227.

\bibitem{ZhangZhu2020} (MR4046741) [10.1016/j.cnsns.2019.105132]
    \newblock G. Zhang, R. Zhu,
    \newblock Runge--Kutta convolution quadrature methods with convergence and stability analysis for nonlinear singular fractional integro-differential equations,
    \newblock \emph{Commun. Nonlinear Sci. Numer. Simul.}, \textbf{84} (2020), 105132.

\bibitem{Zou2019} (MR4003758) [10.1007/s11075-018-0613-0]
     \newblock G. Zou,
     \newblock Numerical solutions to time-fractional stochastic partial differential equations,
     \newblock \emph{Numer. Algorithms}, \textbf{82} (2019), 553--571.

\end{thebibliography}
\end{document}